\documentclass[11pt]{article}
  \usepackage[letterpaper, left=1in, right=1in, top=1in, bottom=1in]{geometry}
  
  \usepackage[affil-it]{authblk}  % for associating author affiliations
  
  \usepackage{palatino}  % font face
  \usepackage{setspace}  % line spacing
  \usepackage{color}  % font color
  \usepackage[dvipsnames]{xcolor}  % adds more color options

	\usepackage[labelfont=bf]{caption}	% caption style

  \usepackage{enumitem}  % control layout of itemize, enumerate, description

  \usepackage{amsmath}  % essential math package
  \usepackage{amssymb}  % math symbols
  \usepackage{amsthm}  % creating proof and theorem environments
  \usepackage{bm}  % boldface letters
  \usepackage{bbm}  % supports \mathbb-style digits, e.g. \mathbbm{1} indicating an all-one vector
  \allowdisplaybreaks[4]  % breaking multiline equations
	\usepackage{relsize}  % can increase the size of math symbols

\usepackage[ruled,vlined]{algorithm2e}
\usepackage{algorithmic}
  
  \DeclareMathOperator*{\minimize}{minimize}
  
  \DeclareMathOperator*{\argmin}{argmin}
  \newcommand{\st}{\mathrm{subject\;to}}

	\theoremstyle{definition}

  \usepackage{pifont}  % dingbats and symbols

	\usepackage{longtable}  % for tables that span more than one page
	\usepackage{multirow}  % for columns spanning multiple rows
	\usepackage{array}  % for column specifications
	\usepackage{booktabs}  % providing nicer tables, using \toprule, \midrule, etc.
  \usepackage{graphicx}  % enhanced support for graphics, allows e.g. \in­clude­graph­ics
  \usepackage{float}  % better positioning of figures, e.g. by using [h!]
  \usepackage{subfig}  % creating subfigures with \subfloat
	\usepackage{hyperref}
	\hypersetup{
		pdfstartview = {FitV},  % fits the width of the page to the window
		colorlinks = true,    % false: boxed links; true: colored links
		linkcolor = NavyBlue,  % color of internal links
		citecolor = NavyBlue,  % color of links to bibliography
		filecolor = NavyBlue,  % color of file links
		urlcolor = NavyBlue  % color of external links
	}  

  \usepackage[round]{natbib}  % for natbib, use bibtex as backend
\usepackage{threeparttable}

\title{BO4IO: A Bayesian optimization approach to inverse optimization with uncertainty quantification}
\author[1]{Yen-An Lu}
\author[1]{Wei-Shou Hu}
\author[2]{Joel A. Paulson \thanks{Corresponding author (paulson.82@osu.edu)}}
\author[1]{Qi Zhang \thanks{Corresponding author (qizh@umn.edu)}}
\affil[1]{Department of Chemical Engineering and Materials Science, University of Minnesota, Minneapolis, MN 55455, USA}
\affil[2]{Department of Chemical and Biomolecular Engineering, The Ohio State University, Columbus, OH 43210, USA}
\date{}

\begin{document}

\maketitle

\begin{abstract}
This work addresses data-driven inverse optimization (IO), where the goal is to estimate unknown parameters in an optimization model from observed decisions that can be assumed to be optimal or near-optimal solutions to the optimization problem. The IO problem is commonly formulated as a large-scale bilevel program that is notoriously difficult to solve. Deviating from traditional exact solution methods, we propose a derivative-free optimization approach based on Bayesian optimization, which we call BO4IO, to solve general IO problems. We treat the IO loss function as a black box and approximate it with a Gaussian process model. Using the predicted posterior function, an acquisition function is minimized at each iteration to query new candidate solutions and sequentially converge to the optimal parameter estimates. The main advantages of using Bayesian optimization for IO are two-fold: (i) it circumvents the need of complex reformulations of the bilevel program or specialized algorithms and can hence enable computational tractability even when the underlying optimization problem is nonconvex or involves discrete variables, and (ii) it allows approximations of the profile likelihood, which provide uncertainty quantification on the IO parameter estimates. We apply the proposed method to three computational case studies, covering different classes of forward optimization problems ranging from convex nonlinear to nonconvex mixed-integer nonlinear programs. Our extensive computational results demonstrate the efficacy and robustness of BO4IO to accurately estimate unknown model parameters from small and noisy datasets. In addition, the proposed profile likelihood analysis has proven to be effective in providing good approximations of the confidence intervals on the parameter estimates and assessing the identifiability of the unknown parameters.
\end{abstract}

\section{Introduction}
% General motivation for inverse optimization and its applications

Inverse optimization (IO) is an emerging method for learning unknown decision-making processes \citep{chan2023inverse}. Following the principle of optimality \citep{schoemaker1991quest}, the fundamental idea is to use mathematical optimization as a model for decision-making, where decisions can be viewed as the optimal or near-optimal solutions to an underlying optimization problem. Given observed decisions made by an agent, the goal of IO is to learn the unknown optimization model that best represents the agent’s decision-making process. A key advantage of the IO approach is that it can directly consider constraints, which allows us to leverage all the modeling flexibility of mathematical programming, incorporate domain knowledge, and hence obtain inherently interpretable decision-making models \citep{gupta2023efficient}; this makes it distinct from common black-box machine learning methods such as deep learning. 

The concept of IO was first introduced in an inverse shortest-path problem, where travel costs in a network are determined based on a given route taken by a user \citep{burton1992instance}. Since then, IO has been applied to problems across a wide variety of fields, including personalized treatment planning in healthcare \citep{chan2014generalized,ghate2020imputing,ajayi2022objective}, network design in transportation and power systems \citep{bertsimas2015data,zhang2018price,fernandez2021inverse}, risk preference learning in portfolio optimization \citep{li2021inverse,yu2023learning}, and inference of cellular decision making in biological systems \citep{burgard2003optimization,uygun2007investigation,zhao2016mapping}. 

% Overview of the solution approaches, hint there is a lack of general algorithms for different classes of optimization problems, and motivate our black-box approach.

Early IO research focused on the development of solution approaches for various classes of deterministic IO problems (IOPs) where it is assumed that observed decisions are exact globally optimal solutions to the unknown forward optimization problem (FOP). \cite{ahuja2001inverse} proposed a general solution method for estimating the cost coefficients of a linear program. Others have proposed extensions that address conic \citep{iyengar2005inverse,zhang2010inverse}, discrete \citep{schaefer2009inverse,wang2009cutting}, and nonlinear convex FOPs \citep{zhang2010augmented}, as well as constraint estimation \citep{chan2020inverse,ghobadi2021inferring}. In recent years, the focus has shifted toward data-driven IO, where parameters are estimated based on a set of noisy observations collected under different experiment (input) conditions \citep{mohajerin2018data}. As such, data-driven IO can address real-world problems where the observed decisions are rarely perfectly optimal given a hypothesized structure of the FOP. 

A data-driven IOP is typically formulated as a bilevel program, featuring as many lower-level problems as there are observations. Each lower-level problem represents an instance of the FOP with the corresponding model inputs and observed decisions. The upper-level problem aims to estimate model parameters that align the FOP solutions with the observed decisions. IOPs are notoriously difficult to solve as their size grows with the number of observations. Common solution approaches involve single-level reformulations based on replacing the lower-level problems with their optimality conditions \citep{keshavarz2011imputing,aswani2018inverse} or cutting plane methods \citep{wang2009cutting}. Recent works have also considered tailored decomposition methods for solving IOPs with high-dimensional FOPs and large datasets \citep{gupta2022decomposition,gupta2023efficient}. In general, existing exact solution methods for IOPs cannot be applied to all classes of IOPs, and designing such algorithms becomes particularly challenging when the FOP is nonconvex.

% Talk about DFO and the relevant work on bilevel or multilevel problems.
In this work, we depart from the existing literature and take a derivative-free optimization (DFO) approach to solve general IOPs, where we treat the evaluation of the loss function as a black box. DFO is commonly applied to optimization problems where the full derivative information is not easily accessible \citep{rios2013derivative}. In the field of bilevel or multilevel optimization, local or global sampled-based algorithms, e.g. evolutionary algorithms, have long been used as metaheuristic solution methods \citep{talbi2013taxonomy,sinha2017review}. In the DFO framework developed by \citet{beykal2020domino}, a bilevel program is solved as a single-level grey-box optimization problem, where at each sampling point, the lower-level problem is solved to generate input-output information that the grey-box solver can leverage to suggest new solution candidates and iteratively optimize the upper-level objective. In this work, we propose to use Bayesian optimization (BO), a model-based DFO method designed to find the optimum of an expensive-to-evaluate objective function within a relatively small number of iterations \citep{frazier2018tutorial}. BO leverages Bayesian statistics and surrogate modeling to update a prior distribution for the objective function as new function evaluations are performed, and it provides a means of balancing exploration and exploitation in sampling the next points. BO has become very popular in recent years and has proven to be an effective method for, for example, hyperparameter tuning in deep learning \citep{snoek2012practical, shahriari2015taking} and design of experiments \citep{greenhill2020bayesian, frazier2016bayesian}. However, only recently, it has also been used to tackle bilevel programs \citep{kieffer2017bayesian, dogan2023bilevel} as well as multilevel optimization problems that arise in feasibility analysis and robust optimization \citep{kudva2022constrained,kudva2024robust}.

% Summary of our work.
In our proposed BO framework for solving general IOPs, which we call BO4IO, we approximate the loss function of the IOP with a non-parametric probabilistic surrogate using Gaussian process (GP) regression \citep{williams2006gaussian}. It converges to the optimal parameter estimates that provide the minimum decision loss by iteratively selecting candidate solutions through the minimization of the acquisition function. Here, the key advantage of BO4IO is that, at each iteration, it only requires the evaluation of the loss function with the current parameter estimates; this can be achieved by directly solving the FOP for every data point, which circumvents the need for a complex reformulation or special algorithm. Consequently, BO4IO remains applicable even when the FOP is nonlinear and nonconvex, and it can consider both continuous and discrete decision variables. This approach also naturally decomposes the problem in a sense that the FOP can be solved independently for each data point, which allows parallelization and hence enhances the scalability of the algorithm with respect to the number of data points. 

% Motivation on PL
BO4IO also enables us to address another major challenge in IO that is often overlooked, namely the quantification of uncertainty with respect to the estimated parameters. In IO, there are often multiple parameter estimates that lead to the same loss; these ``equally good'' estimates form the so-called inverse-feasible set. The size of the inverse-feasible set can be viewed as a measure of uncertainty in any point estimate drawn from that set; the larger the inverse-feasible set, the less confident one can be about the estimate. In the case of linear programs where the cost vector is unknown, the inverse-feasible set takes the form of a polyhedral cone; \citet{gupta2022decomposition} use this insight to develop an adaptive sampling method that tries to sample points that can help further reduce the size of the inverse-feasible set. However, that method is specific to linear programs and cannot be directly extended to general IOPs. A common approach in parameter estimation is to derive confidence intervals for the point estimates. Asymptotic confidence intervals can be computed using a Fisher information matrix \citep{sobel1982asymptotic}, but it is only exact when the estimated model is linear \citep{joshi2006exploiting}. In contrast, sample-based confidence intervals derived using the profile likelihood function also apply to the nonlinear case \citep{neale1997use,raue2009structural}. We find that in BO4IO, we can use the posterior of the GP surrogate to approximate the profile likelihood, which provides uncertainty quantification and insights into the identifiability of given model parameters. 

To assess the efficacy of the proposed BO approach for various classes of FOPs, we perform computational case studies wherein we learn cellular objectives in flux balance analysis of metabolic networks as well as market requirements in standard and generalized pooling problems. The computational experiments demonstrate that BO4IO can efficiently identify estimates of the unknown model parameters with a relatively small number of iterations. It proves to be especially beneficial in cases where the size of the FOP is large relative to the number of unknown parameters. In addition, using the approximate profile likelihood, we can assess how the identifiability of certain parameters changes with the number and quality of observations. To the best of our knowledge, this is the first DFO method specifically designed to solve IOPs as well as the first approach that provides a direct measure of uncertainty on the IO parameter estimates.

The remainder of this paper is structured as follows. In Section \ref{sec:problem_formulation}, we present a general formulation of an IOP. We provide the necessary background on GPs and BO and introduce our BO4IO algorithm in Section \ref{sec:BO4IO}. In Section \ref{sec:PL}, we introduce the profile likelihood and derive its approximation based on the GP surrogate in BO4IO. Finally, we demonstrate the performance of BO4IO in three numerical case studies in Section \ref{sec:case_studies} and conclude in Section \ref{sec:conclusions}. 

\section{The inverse optimization problem}\label{sec:problem_formulation}
% Definition of FOPs
In IO, we assume that a decision-making process can be modeled as an optimization problem, also referred to as the \textit{forward} optimization problem, formulated in the following general form:
\begin{equation}
\label{eqn:FOP}
\tag{FOP}
\begin{aligned}
  \minimize_{x \in \mathbb{R}^n} \quad & f(x,u;\theta) \\
  \st \quad & g(x,u;\theta) \leq 0,
\end{aligned}
\end{equation}
where $x \in \mathbb{R}^n$ is the vector of decision variables, $u \in \mathbb{R}^m$ denotes the vector of contextual inputs that describe the system conditions, and $f$ and $g$ represent the objective and constraint functions, respectively, that are parameterized by model parameters $\theta \in \Theta$, where $\Theta$ is some compact set in $\mathbb{R}^d$.
% General formulation of IOPs
In a data-driven IO setting, $\theta$ are unknown, but decisions $\{x_i\}_{i\in\mathcal{I}}$ can be collected under varying conditions $\{u_i\}_{i\in\mathcal{I}}$, where $\mathcal{I}$ denotes the set of observations. The observations are assumed to be noisy, which could be due to, for example, measurement errors, suboptimal decisions, or model mismatch. The goal of IO is to estimate $\theta$ such that the solutions of \eqref{eqn:FOP} best fit the observed decisions, which gives rise to the following bilevel optimization problem:
\begin{equation}
\label{eqn:IOP}
\tag{IOP}
\begin{aligned}
  \minimize_{\hat{\theta} \in \Theta, \, \hat{x}} \quad & l(\hat{\theta}):= \sum_{i \in \mathcal{I}} (x_{i} - \hat{x}_{i}(\hat{\theta}))^{\top} W (x_{i} - \hat{x}_{i}(\hat{\theta})) \\
  \st \quad & \hat{x}_{i}(\hat{\theta}) \in \argmin_{\tilde{x} \in \mathbb{R}^n} \left\lbrace f(\tilde{x},u_i;\hat{\theta}): g(\tilde{x},u_i;\hat{\theta}) \leq 0 \right\rbrace \quad \forall \, i \in \mathcal{I},
\end{aligned}
\end{equation}
where the objective of the upper-level problem is to obtain parameter estimates, denoted by $\hat{\theta}$, such that the (empirical) decision loss $l : \Theta \to \mathbb{R}_+$, a function that quantifies the agreement between the model predictions ($\hat{x}_{i}$) and observations ($x_{i}$), is minimized. For simplicity, we focus on the case where $l$ is defined as a weighted sum of squared residuals, where $W \in \mathbb{S}^{n}_+$ denotes the positive definite matrix of weighting factors; this loss function can be easily replaced with a more general likelihood or posterior loss function (see, e.g., \cite[Chapter 5]{Goodfellow-et-al-2016}). In the $|\mathcal{I}|$ lower-level problems of \eqref{eqn:IOP}, one for each observation, $\hat{x}_i$ is constrained to be an optimal solution to the corresponding optimization problem parameterized by $\hat{\theta}$; as such, $\hat{x}_i$ is a function of $\hat{\theta}$.

\section{The BO4IO framework}\label{sec:BO4IO}
In this section, we provide a brief overview of derivative-free optimization and Gaussian process regression, and show how they are applied to solve \eqref{eqn:IOP} in the proposed Bayesian optimization for inverse optimization (BO4IO) framework.

\subsection{Derivative-free inverse optimization}
We propose to solve general IOPs using a DFO approach where we treat \eqref{eqn:IOP} as the following black-box optimization problem:
\begin{equation*}
\minimize_{\hat{\theta} \in \Theta} \quad l(\hat{\theta}).
\end{equation*}
DFO methods aim to sequentially optimize the target objective within a designated search space, where the query points are suggested based on past function evaluations. At each iteration, the objective function, which in our case is the loss function $l$, is evaluated for a given $\hat{\theta}$ by solving the $|\mathcal{I}|$ individual FOPs to global optimality. Specifically, for every $i \in \mathcal{I}$, we solve \eqref{eqn:FOP} with $u = u_i$ and $\theta = \hat{\theta}$ to obtain $\hat{x}_i(\hat{\theta})$ with which we can then compute $l(\hat{\theta}) = \sum_{i \in \mathcal{I}} (x_{i} - \hat{x}_{i}(\hat{\theta}))^{\top} W (x_{i} - \hat{x}_{i}(\hat{\theta}))$. Note that each FOP can be solved independently; hence, they can be solved in parallel to further reduce the computation time.

Even though parallel computing can be applied, each evaluation of the loss function can still be expensive if the underlying FOP is difficult to solve; hence, it is important to select a DFO algorithm that requires as few evaluations as possible. DFO methods can be roughly categorized into direct and model-based methods \citep{rios2013derivative}. While direct algorithms explicitly determine a search direction from previous function evaluations, model-based methods use these function evaluations to construct a surrogate model that approximates the target function and determine the next sampling point using that surrogate \citep{jones1998efficient,huyer2008snobfit}. %By providing an accurate approximation of functional space, the model-based approach can converge to nearly global solutions with a restricted number of iterations. 
In this work, we choose Bayesian optimization (BO) to be our preferred model-based DFO method \citep{mockus1994application,frazier2018tutorial}. In BO, a probabilistic surrogate model is constructed where Bayes' rule is used to compute the surrogate posterior distribution based on the selected prior and past function evaluations. BO uses the posterior to quantify the uncertainty in the estimate and forms a utility (acquisition) function that balances the exploration-exploitation trade-off in the search process. A wide variety of probabilistic surrogate models have been used in BO, e.g., Gaussian processes \citep{mockus1994application}, random forests \citep{hutter2011sequential}, tree-structured Parzen estimators \citep{bergstra2011algorithms}, and deep neural networks \citep{snoek2015scalable}. In this work, we use GP as the probabilistic surrogate due to its nonparametric feature and the analytical expression of the posterior; however, it should be noted that the proposed framework can readily consider other surrogate models.

\subsection{Gaussian process regression}
Based on the standard GP detailed in \citet{williams2006gaussian}, we build a GP model that approximates the loss function $l: \Theta \to \mathbb{R}_+$. The approximation is obtained by assuming that the function values $l(\hat{\theta})$ at different $\hat{\theta}$ values are random variables of which any finite subset forms a joint Gaussian distribution. A GP model can thus be interpreted as an infinite collection of these random variables, describing a distribution over the function space.

We let $\hat{l}(\hat{\theta})\sim \mathcal{GP}(\mu(\hat{\theta}), \kappa(\hat{\theta}, \hat{\theta}'))$ denote the GP surrogate of the true loss function $l$. The GP model is fully specified by a prior mean function $\mu(\hat{\theta})$ and a covariance (kernel) function $\kappa(\hat{\theta}, \hat{\theta}')$ that is the covariance of the function values $\hat{l}(\hat{\theta})$ and $\hat{l}(\hat{\theta}')$ for any $\hat{\theta}, \, \hat{\theta}' \in \Theta$. Without loss of generality, we assume that $\mu(\hat{\theta}) = 0$, which can be achieved by normalizing the output data when fitting the surrogate. The kernel function encodes the smoothness and rate of change in the target function (see, e.g., \citet[Chapter 4]{williams2006gaussian}). We focus on the stationary covariance function from the Matern class, which is a popular kernel choice in many GP applications.  Let $\mathcal{D}_t = \{(\hat{\theta}_i, l(\hat{\theta}_i)\}^t_{i=1}$ denote a set of $t$ past evaluations. Conditioned on the evaluations in $\mathcal{D}_t$, the predicted posterior distribution of the function $\hat{l}$ at a future input  $\hat{\theta}_{t+1}$ remains Gaussian with the following posterior mean $\mu_t(\hat{\theta}_{t+1})$ and covariance $\sigma^2_t(\hat{\theta}_{t+1})$:
\begin{equation*}
\begin{aligned}
  & \mu_t(\hat{\theta}_{t+1}) = \boldsymbol{\kappa}_t^{\top}(\hat{\theta}_{t+1})\boldsymbol{K}_t^{-1}\boldsymbol{l}_t\\
  & \sigma^2_t(\hat{\theta}_{t+1}) = \kappa(\hat{\theta}_{t+1},\hat{\theta}_{t+1})-\boldsymbol{\kappa}_t^\top(\hat{\theta}_{t+1})\boldsymbol{K}_t^{-1}\boldsymbol{\kappa}_t(\hat{\theta}_{t+1}),
\end{aligned}
\end{equation*}
where $\boldsymbol{l}_t=[l(\hat{\theta}_1),...,l(\hat{\theta}_t)]^\top$ is the vector of $t$ past loss function evaluations, $\boldsymbol{K}_t \in \mathbb{R}^{t\times t}$  is the covariance matrix with entries $[\boldsymbol{K}_t]_{i,j}=\kappa(\hat{\theta}_i,\hat{\theta}_j)$ for all $i,j \in \{1,...,t\}$, and $\boldsymbol{\kappa}(\hat{\theta}_{t+1})$ is the vector of covariance values between the new input $\hat{\theta}_{t+1}$ and the evaluated inputs $\{\hat{\theta}_i\}_{i=1}^t$. In practice, the kernel function contains unknown hyperparameters, such as smoothness and length-scale parameters, which are calibrated to the dataset $\mathcal{D}_t$ using maximum likelihood estimation. 

\subsection{The BO4IO algorithm}
A pseudocode for the proposed BO4IO algorithm is shown in Algorithm \ref{alg:BO4IO}, which aims to minimize the decision loss $l(\hat{\theta})$ of \eqref{eqn:IOP} over the parameter domain $\hat{\theta}\in\Theta$ by sequentially querying $\hat{\theta}$ based on the predicted posterior distribution of the GP surrogate. At the $t$th iteration, the GP model is updated using the past $t$ evaluations ($\mathcal{D}_t$). An acquisition function that encodes the exploration-exploitation trade-off is constructed based on the predicted GP posterior and is optimized to query the next sampling point $\hat{\theta}_{t+1}$. The true loss function value $l(\hat{\theta}_{t+1})$ is then evaluated by solving $|\mathcal{I}|$ FOPs. The dataset $\mathcal{D}_t$ is then appended with the new evaluation $(\hat{\theta}_{t+1}, l(\hat{\theta}_{t+1}))$ to create a concatenated dataset $\mathcal{D}_{t+1}$. The new dataset  $\mathcal{D}_{t+1}$ is then used to update the GP model, and the same process is repeated until the maximum number of iterations is reached (or some other termination criterion is met). 

A wide variety of acquisition functions have been developed in BO; here, we use the lower confidence bound (LCB) acquisition function \citep{auer2002using} due to its simplicity and nice theoretical properties. The acquisition function is minimized to obtain the next sampling point, i.e.
\begin{equation}
\label{eqn:lcb}
\begin{aligned}
  \hat{\theta}_{t+1} \in \argmin_{\hat{\theta} \in\Theta} \ \mu_t(\hat{\theta}) - \beta^{1/2}_t\sigma_t(\hat{\theta}),
\end{aligned}
\end{equation}
where $\beta_t > 0$ is an exploration factor. When $\beta_t$ is defined as an iteration-dependent parameter, some theoretical guarantees can be derived for the convergence to a global optimal solution \citep{srinivas2009gaussian}; however, the resulting search process may become very slow with the adaptive approach. Instead, one often considers a constant $\beta_t$ value that bounds the failure probability per iteration. For example, a value of $\beta_t = 4$ has demonstrated good results in practical applications \citep{berkenkamp2019no}. 

\begin{algorithm}
\caption{BO4IO: Bayesian optimization for inverse optimization}
\label{alg:BO4IO}
\textbf{Input:} \begin{tabular}[t]{l}
 parameter domain $\Theta$;\\
 kernel for GP prior $\kappa(\hat{\theta},\hat{\theta}')$;\\
 exploration parameter $\beta$;\\
 set of observations $\{(u_i,x_i)\}_{i=1}^{|\mathcal{I}|}$;\\
 FOPs with inputs $\{u_i\}_{i=1}^{|\mathcal{I}|}$;\\
 initial evaluation dataset $\mathcal{D}_0$;\\
 maximum number of iterations $T$.
\end{tabular}\\
\begin{algorithmic}[1]
\STATE Generate GP model using $\mathcal{D}_0$
\FOR{$t=1,2,...,T$}
    \STATE $\hat{\theta}_{t+1} \in \argmin_{\hat{\theta}\in\Theta} \mu_t(\hat{\theta}) - \beta \sigma_t(\hat{\theta})$, see \eqref{eqn:lcb}
    \STATE Solve the $\mathcal{|I|}$ FOPs with $\hat{\theta}_{t+1}$ and compute $l(\hat{\theta}_{t+1})$
    \STATE $\mathcal{D}_{t+1} \leftarrow \mathcal{D}_t\cup\{(\hat{\theta}_{t+1}, l(\hat{\theta}_{t+1}))\}$
    \STATE Update GP model using $\mathcal{D}_{t+1}$
\ENDFOR
\RETURN{optimal solution $\hat{\theta}^*_T \in \argmin_{\hat{\theta}_t \in \mathcal{D}_T}\{l(\hat{\theta_t})\}_{t=1}^T$} 
\end{algorithmic}
\end{algorithm}
% \subsection{The algorithm}

\section{Uncertainty quantification using profile likelihood}\label{sec:PL}
In this section, we propose a method to quantify the uncertainty in the parameter estimates from BO4IO using the profile likelihood (PL) method. We first review uncertainty quantification in general parameter estimation problems and discuss how the confidence interval (CI) of each parameter estimate can be derived using the PL method. We then propose approximations of the PLs and CIs based on the GP surrogate in BO4IO, which can be used to derive optimistic and worst-case CIs of parameter estimates and provide insights into the identifiability of the unknown parameters in IO.

\subsection{General approach for parameter estimation and uncertainty quantification}
It is widely known that the weighted sum of the residuals function $l(\theta)$ defined in \eqref{eqn:IOP} is directly related to the likelihood function under the assumption of independent and identically distributed (i.i.d.) Gaussian observation noise. In particular, if $x_i = \hat{x}_i(\hat{\theta}) + \epsilon_i$ with $\epsilon_i \sim \mathcal{N}( 0, \sigma_i^2 )$, then
\begin{equation*}
    l(\hat{\theta}) = c - 2\log( L(\hat{\theta}) ),
\end{equation*}
where $L(\hat{\theta}) = p(\{ x_i \}_{i \in \mathcal{I}} | \{ u_i \}_{i \in \mathcal{I}}, \hat{\theta})$ is the likelihood function, $c$ is a constant, and $l(\hat{\theta})$ is given in \eqref{eqn:IOP} with $W = \text{diag}( \sigma_1, \ldots, \sigma_{|\mathcal{I}|} )$. Let $\hat{\theta}^* \in \argmin_{\hat{\theta} \in \Theta} l(\hat{\theta})$ denote the maximum likelihood estimate (MLE) of the unknown parameters $\theta$. Finite sample confidence intervals, which signify that the true value is within a set with at least some probability $1-\alpha$, can be derived as follows
\begin{align}\label{eqn:CI}
    \text{CI}_{\hat{\theta}^*} = \{ \hat{\theta} \in \Theta : l(\hat{\theta}) - l(\hat{\theta}^*) < \Delta_\alpha \},
\end{align}
where $\Delta_\alpha$ denotes the $\alpha$-level quantile of the standard $\chi^2$ distribution with degrees of freedom equal to 1 for pointwise and $d$ for simultaneous CIs, respectively. 

An important challenge in many parameter estimation problems is determining the \textit{identifiability} of the parameters to be estimated in the model \citep{bellman1970structural,cobelli1980parameter}. We can think of identifiability in terms of the size of the CI. If the CI is small, we have found a narrow range in which the true parameter exists and so we can say that the model is fully identifiable (and we can consequently trust the found estimates $\hat{\theta}^*$). On the other hand, if the CI is extended in any direction, we are unable to pinpoint an accurate estimate of the parameter; hence, the parameter is non-identifiable. As discussed in \cite{wieland2021structural}, there are two types of non-identifiabilities that can emerge in a model. The first is \textit{structural non-identifiability}, which is the result of an effectively redundant model parametrization. Structurally non-identifiable parameters result in non-unique MLE estimates $\hat{\theta}^*$ such that they cannot be resolved by gathering more data or improving the quality of the existing data. The second is \textit{practical non-identifiability}, which refers to a parameter that is structurally identifiable (in the sense that it has a unique MLE estimate) but $l(\hat{\theta}) - l(\hat{\theta}^*) < \Delta_\alpha$ remains below the threshold $\Delta_\alpha$ for a large portion of $\Theta$. One may be able to make practically non-identifiable parameters identifiable by collecting additional and/or higher-quality data to reduce the CI. One can use model-based design of experiments (DOE) methods \citep{balsa2008computational,bandara2009optimal,kreutz2009systems} to systematically select new experimental trials that are most informative for reducing the CI. We plan to study the DOE problem for \eqref{eqn:IOP} in more detail in future work. 

\subsection{The profile likelihood method}
The PL method is an effective strategy for differentiating between the different identifiability categories, even in high-dimensional spaces \citep{raue2009structural}. The approach assesses the identifiability of a parameter of interest $\hat{\theta}_k$ by profiling the maximum likelihood (i.e. minimum loss). The profile likelihood $\text{PL}(\hat{\theta}_k)$ for a given value of $\hat{\theta}_k$ is defined as
\begin{gather}\label{eqn:PL}
    \text{PL}(\hat{\theta}_k) = \min_{\hat{\theta}_{k'\neq k}} \ l(\hat{\theta}).
\end{gather}
Note that $\text{PL}(\hat{\theta}_k)$ is obtained by fixing $\hat{\theta}_k$ while re-optimizing the rest of the parameters. Based on the definition given in \eqref{eqn:CI}, the CI of each parameter estimate $\theta^*_k$ can be computed as follows:
\begin{gather}
    \text{CI}_{\hat{\theta}^*_k} = \{ \hat{\theta}_k \mid \text{PL}(\hat{\theta}_k) - l(\hat{\theta}^*)\leq  \Delta_\alpha \},
\end{gather}
where $l(\hat{\theta}^*)$ is the decision loss with the optimal parameter estimates of \eqref{eqn:IOP}. It should be noted that the re-optimization during the PL evaluation provides useful information about the nonlinear relationships between the parameters, which can be further leveraged in model reduction \citep{maiwald2016driving}.

\subsection{Approximation of profile likelihood and confidence intervals in BO4IO}
\begin{figure}[hp]
\centering
\includegraphics[scale=0.55]{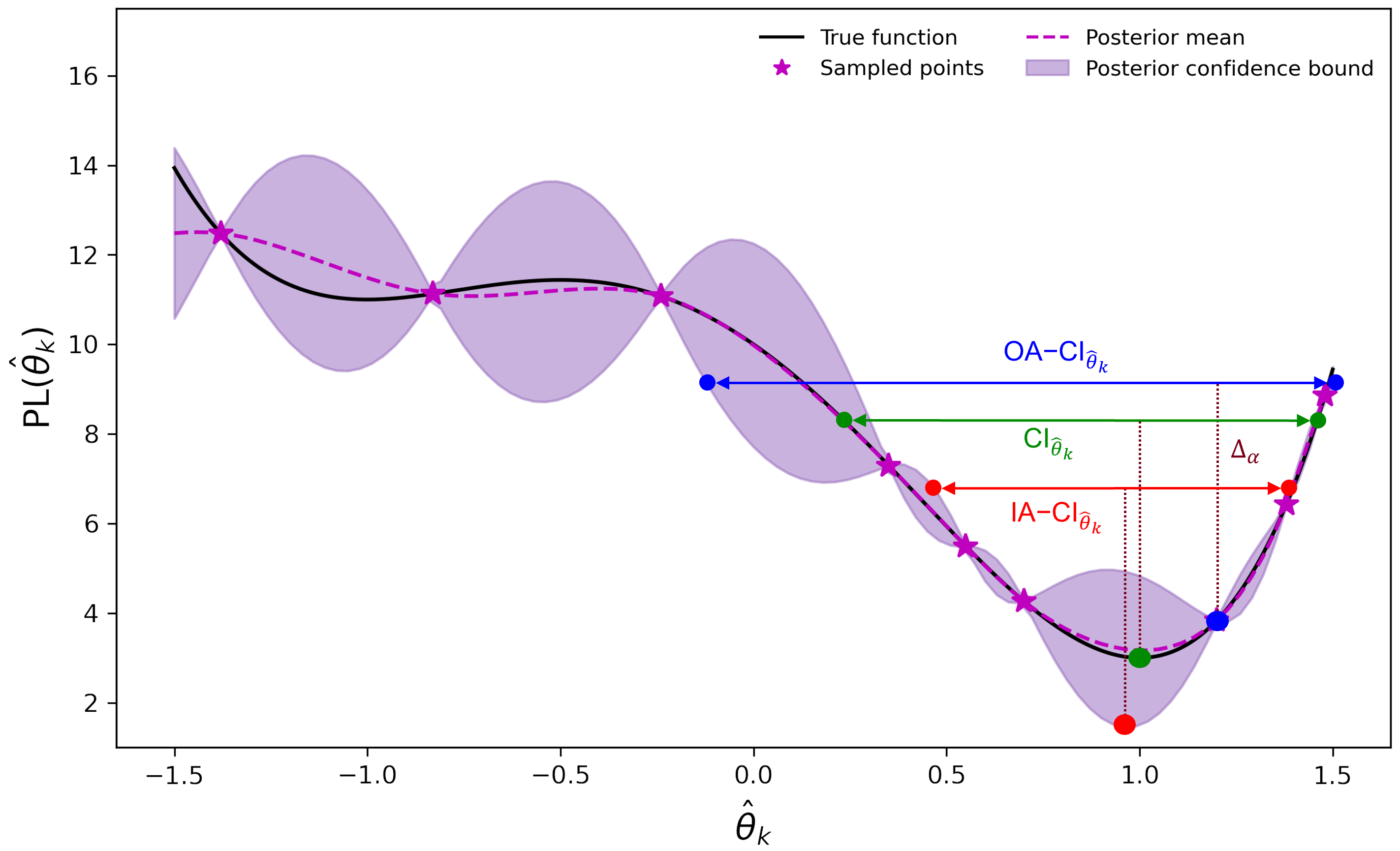}
\caption{Illustrative example of outer-approximate (OA-CI, blue), original (CI, green), and inner-approximate (IA-CI, red) confidence intervals of a parameter of interest ($\hat{\theta}_k$) based on the proposed PL analysis.}
\label{fig:PL_exmple}
\end{figure}

Similar to \eqref{eqn:IOP}, solving \eqref{eqn:PL} exactly is challenging for general FOPs; as a result, we do not have access to the exact PL values and CIs of the parameter estimates. Instead, we approximate them using the GP posterior obtained in BO4IO. We first define the upper and lower confidence bounds on the true loss function $l(\hat{\theta})$ at the $t$th iteration as follows:
\begin{align*}
    & \hat{l}^\text{UCB}_t(\hat{\theta}) = \mu_t(\hat{\theta}) + \rho_t^{1/2} \sigma_t(\hat{\theta}) \\
    & \hat{l}^\text{LCB}_t(\hat{\theta}) = \mu_t(\hat{\theta}) - \rho_t^{1/2} \sigma_t(\hat{\theta}),
\end{align*}
where $\rho_t$ is a parameter denoting the confidence levels. The confidence bounds on the PL of a parameter of interest $\hat{\theta}_k$ at the $t$th iteration can be further derived as
\begin{align}
    & \widehat{\text{PL}}^\text{UCB}_t(\hat{\theta}_k) = \min_{\hat{\theta}_{k'\neq k}} \ \hat{l}^\text{UCB}_t(\hat{\theta}) \label{eqn:PL-UCB} \\
    & \widehat{\text{PL}}^\text{LCB}_t(\hat{\theta}_k) = \min_{\hat{\theta}_{k'\neq k}} \ \hat{l}^\text{LCB}_t(\hat{\theta}), \label{eqn:PL-LCB}
\end{align}
where $\widehat{\text{PL}}^\text{UCB}_t(\hat{\theta}_k)$ and $\widehat{\text{PL}}^\text{LCB}_t(\hat{\theta}_k)$ denote the upper and lower confidence bounds on the PL, respectively. In addition, the optimistic estimate of the minimum decision loss can be determined by $\hat{l}^\text{*LCB}_t = \min l^\text{LCB}_t(\hat{\theta})$, whereas a pessimistic estimate of the minimum decision loss is the current best-found loss $l^*_t = \min\{l(\hat{\theta}_i),...,l(\hat{\theta}_t)\}$. The outer-approximation (OA) and inner-approximation (IA) of the CI on each best parameter estimate $\hat{\theta}^*_k$ at the $t$th iteration can hence be written as 
\begin{align}
    & \text{OA-CI}_{t,\hat{\theta}^*_k} = \{ \hat{\theta}_k | \widehat{\text{PL}}^\text{LCB}_t(\hat{\theta}_k) - l^*_t\leq  \Delta_\alpha \} \label{eqn:OA-CI} \\
    & \text{IA-CI}_{t,\hat{\theta}^*_k}= \{ \hat{\theta}_k | \widehat{\text{PL}}^\text{UCB}_t(\hat{\theta}_k) - \hat{l}^\text{*LCB}_t\leq  \Delta_\alpha \}.
\end{align}
Since $l^*_t \geq l^* \geq \hat{l}_t^{*\text{LCB}}$ and $\widehat{\text{PL}}^\text{UCB}_t(\hat{\theta}_k) \geq \text{PL}(\hat{\theta}_k) \geq \widehat{\text{PL}}^\text{LCB}_t(\hat{\theta}_k)$, we have $\text{OA-CI}_{t,\hat{\theta}^*_k} \subseteq \text{CI}_{\hat{\theta}^*_k}\subseteq \text{IA-CI}_{t,\hat{\theta}^*_k}$. In other words, $\text{OA-CI}_{t,\hat{\theta}^*_k}$ provide a worst-case estimate of the confidence region, whereas $\text{IA-CI}_{t,\hat{\theta}^*_k}$ is an optimistic estimate. An illustrative example of the approximate and true CIs are shown in Figure \ref{fig:PL_exmple}. Notably, we are primarily interested in the worst-case estimate of the CI; hence, we only report $\text{OA-CI}_{t,\hat{\theta}^*_k}$ in our case studies. The pseudocode for the algorithm used to determine $\text{OA-CI}_{t,\hat{\theta}^*_k}$ is shown in Algorithm \ref{alg:OACI}.

\begin{algorithm}
\caption{Outer-approximation of confidence intervals in BO4IO}
\label{alg:OACI}
\textbf{Input:} \begin{tabular}[ht!]{l}
 parameter of interest, $\hat{\theta}_k$;\\
 posterior mean and covariance function of GP $\mu_t(\hat{\theta})$ and $\sigma_t^2 (\hat{\theta})$;\\
 confidence parameter $\rho$;\\
 the $\alpha$-quantile of $\chi^2$ distribution $\Delta_\alpha$;\\
 the best decision loss of $t$ past evaluations $l^*_t$;\\
 the step size for $\hat{\theta}_k$ $\delta_{\hat{\theta}_k}$;\\
 the sampled PL range $[\hat{\theta}_k^\text{LB}, \hat{\theta}_k^\text{UB}]$ 
\end{tabular}\\
% \textbf{Output:} 
% \begin{tabular}[t]{l}
% $\hat{\text{CI}}_{t,\theta_k} = \{ \theta_k | \hat{l}_t^\text{PL}(\theta_k) - \hat{l}_t^*\leq  \Delta_\alpha \}$
% \end{tabular}
\begin{algorithmic}[1]
\STATE $\mathcal{S} \leftarrow  \emptyset$ and $\hat{\theta}_k \leftarrow  \hat{\theta}_k^\text{LB}$
\WHILE{$\hat{\theta}_k \leq \hat{\theta}_k^\text{UB}$}
    \STATE $\widehat{\text{PL}}^\text{LCB}_t(\hat{\theta}_k) \leftarrow \min_{\hat{\theta}_{k'\neq k}} \mu_t(\hat{\theta}) - \rho \sigma_t(\hat{\theta})$, see \eqref{eqn:PL-LCB}
    \STATE $\mathcal{S} \leftarrow \mathcal{S}\cup\{(\hat{\theta}_k, \, \widehat{\text{PL}}^\text{LCB}_t(\hat{\theta}_k))\}$
    \STATE $\hat{\theta}_k=\hat{\theta}_k+\delta_{\hat{\theta}_k}$
\ENDWHILE
\STATE $\text{OA-CI}_{t,\hat{\theta}^*_k} \leftarrow \{ \hat{\theta}_k | \widehat{\text{PL}}^\text{LCB}_t(\hat{\theta}_k) - \hat{l}_t^*\leq  \Delta_\alpha \}$ based on $\mathcal{S}$
\RETURN{$\text{OA-CI}_{t,\hat{\theta}^*_k}$}
\end{algorithmic}
\end{algorithm}
% \newpage
\section{Computational case studies}\label{sec:case_studies}
We apply the proposed BO4IO algorithm to three computational case studies covering various classes of FOPs. In the first case study, we estimate the optimal linear combinations of multiple cellular objectives in the flux balance analysis, which is formulated as a convex nonlinear program (NLP). In the second case study, we learn the product demands in the standard pooling problem, which is formulated as a nonconvex NLP. Finally, we extend the second study to a generalized pooling problem that contains mixed-integer decisions and estimate the product quality constraints. Using synthetic data, we evaluate the performance of the algorithm and apply the approximate PL method to assess the confidence intervals and identifiability of the IO parameter estimates. For each specific case, 10 random instances were considered to provide reliable computational statistics. The algorithm was implemented in Python 3.8 using Pyomo \citep{bynum2021pyomo} as the optimization modeling environment and Gurobi v10.0.1 \citep{gurobi} as the optimization solver. The BO framework was built under BoTorch \citep{balandat2020botorch} where GPyTorch \citep{gardner2018gpytorch} was used to construct the GP models. All codes were run on the Agate cluster of the Minnesota Supercomputing Institute (MSI) using a node allocated with 32 GB of RAM and 64 cores. The implementation of all case studies can be found in our GitHub repository at \href{https://github.com/ddolab/BO4IO}{https://github.com/ddolab/BO4IO}.

\subsection{Case study 1: Learning cellular objectives in flux balance analysis}
Flux balance analysis (FBA) is a method for simulating the flux distributions in the metabolic network of a cell at steady state \citep{orth2010flux}. Metabolic flux is one of the critical determinants of cellular status under different physiological and environmental conditions; the FBA method serves as a powerful tool for the optimization of biosynthesis in microbial industries \citep{pardelha2012flux, choon2014differential} and the identification of drug targets in human diseases \citep{lee2009comparative,nilsson2017genome}. The original FBA problem is formulated as a linear program by assuming that the cells try to optimize a cellular objective subject to the reaction stoichiometry of cell metabolism. An extension with $l_1$ or $l_2$ regularization is often used to avoid degeneracy in the FBA solutions \citep{bonarius1996metabolic,lewis2010omic}. Moreover, a variation that considers multi-objective optimization via a weighted combination of common objectives has proven to better describe the fluxes in some organisms \citep{lee2004identification, schuetz2012multidimensional, garcia2012predictive}; however, the weights are often arbitrarily assigned \citep{nagrath2010soft}. The goal of this IOP is to learn the weight vector in multi-objective FBA from the synthetic datasets. The FBA problem is formulated as the following quadratic FOP with $l_2$-regularization:
\begin{subequations}
\label{eqn:Convex-MOO}
% \tag{Convex-MOO}
\begin{align}
      \minimize_{\boldsymbol{v}} \quad & \sum_{k\in \mathcal{R}^\text{obj}}\theta_k v_k + \lambda \sum_{k \in \mathcal{R}} v_k^2 \label{eqn:Convex-MOO-a}\\
 \mathrm {subject \; to} \quad &\sum_{k \in \mathcal{R}} s_{jk} v_k = 0\quad \forall j \in \mathcal{M}\label{eqn:Convex-MOO-b}\\
\quad & L_k \leq v_k \leq U_k \quad \forall k \in \mathcal{R}, \label{eqn:Convex-MOO-c}
\end{align}
\end{subequations}
where $\mathcal{M}$ and $\mathcal{R}$ denote the sets of metabolites and metabolic reactions, respectively. Here, the reaction fluxes are represented by $v$. The first term of the objective function \eqref{eqn:Convex-MOO-a} is a linear combination of common cellular objectives in the set $\mathcal{R}^\text{obj}$ with unknown weight parameters $\theta$. The second term of \eqref{eqn:Convex-MOO-a} is the $l_2$-regularization scaled by a factor $\lambda$. Constraints \eqref{eqn:Convex-MOO-b} represent the material balances for all metabolites, where $s_{jk}$ denotes the stoichiometric coefficient of metabolite $j$ in reaction $k$. Constraints \eqref{eqn:Convex-MOO-c} set lower and upper bounds, respectively denoted by $L$ and $U$, on the metabolic fluxes. 

We consider the core \textit{E. coli} model (ID: e\textunderscore coli\textunderscore core) from the BiGG database \citep{orth2010reconstruction,king2016bigg}, where the metabolic network specification ($|\mathcal{M}|$, $|\mathcal{R}|$) is $(72,95)$. The true cellular objective is assumed to be an unknown linear combination of up to five common cellular objectives in $\mathcal{R}^\text{obj}$, including maximization of biomass formation, maximization of ATP generation, minimization of $\text{CO}_\text{2}$, minimization of nutrient uptake, and minimization of redox potential \citep{savinell1992network,schuetz2007systematic,garcia2014comparison} with an unknown weight vector $\theta$. In each random instance of the IOP, we first generate the ground-truth parameters $\theta$ from a Dirichlet distribution, where the weight vector is scaled by assuming the summation of all weights to be one. The same equality constraint is considered when minimizing the acquisition function; hence, the number of unknown parameters is $d = |\mathcal{R}^\text{obj}|-1$. Next, using the same $\theta$ vector, we collect the optimal flux vectors $v^*_i$ for each observation $i$ by solving the FOP under varying input conditions $\{u_i\}_{i\in\mathcal{I}} = \{(L_i, U_i)\}_{i\in\mathcal{I}}$, where the flux bound vectors $L_i, \, U_i \in \mathbb{R}^{|\mathcal{R}|}$ are sampled from $L_{ik}, U_{ik} \sim \mathcal{U}(10,100)$ for each $k \in\mathcal{R}$. We then generate the noisy data by first standardizing each optimal flux $v^*_{ik}$ across all observations and perturb the standardized values such that $\bar{v}_{ik} = \bar{v}^*_{ik} + \epsilon_{ik}$ for all $i\in \mathcal{I}$ and $k\in \mathcal{R}$, where $\bar{v}^*_{ik}$ denotes the standardized optimal flux for reaction $k$ and $\epsilon_{ik} \sim \mathcal{N}(0,\sigma^2)$. 

\subsubsection{Computational performance under varying dimensionality}
\begin{figure}[h]
\centering
\includegraphics[scale=0.54]{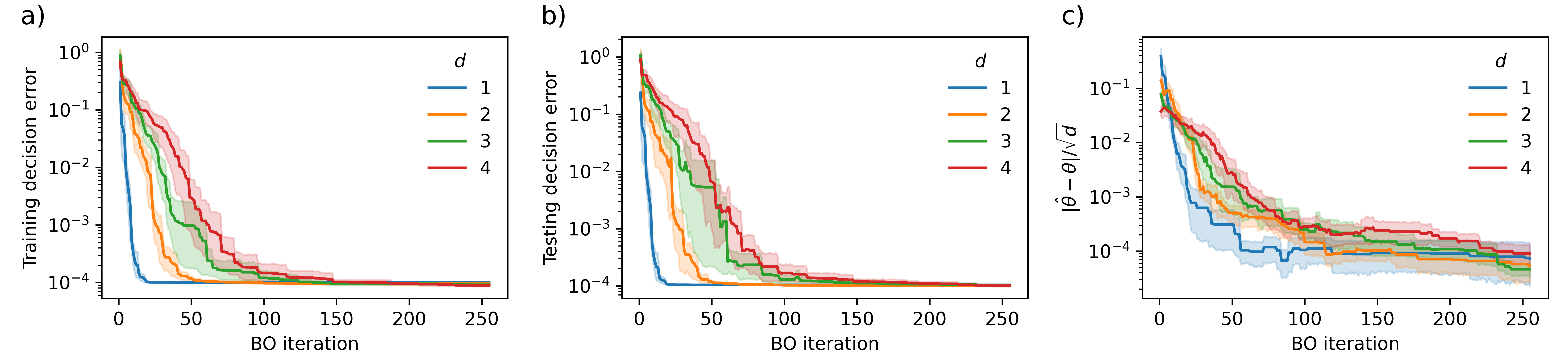}
\caption{Effect of the dimensionality ($d$) of $\theta$ on the accuracy of the estimated FOPs. Convergence analysis of (a) training error, (b) testing error, and (c) parameter error with varying $d$. Training and testing errors refer to the average standardized decision error defined as $\sum_{i\in\mathcal{I}}\sum_{k\in{\mathcal{R}}}(\bar{v}_{ik}-\hat{v}_{ik})^2/|\mathcal{I}|/|\mathcal{R}|$ and calculated based on the training and testing datasets, respectively, whereas the parameter error denotes the difference between the ground-truth ($\theta$) and the estimated ($\hat{\theta}^*$) values. Here, the solid lines and shaded areas respectively denote the medians and confidence intervals of the corresponding error across the 10 random instances. The synthetic dataset of each random instance is generated under the setting of $|\mathcal{I}| = 50$ and $\sigma=0.01$.}
\label{fig:FBA_diff_dim}
\end{figure}

We run 10 random instances with 50 training data points, i.e. $|\mathcal{I}|=50$, under varying dimensionality of unknown parameters ($d=1 \sim 4$) at a low noise level of $\sigma=0.01$. Separate testing datasets with the same $|\mathcal{I}|$ number of data points are generated to test the model prediction. The convergence profiles for the training and testing decision errors as well as the parameter error over 250 BO iterations are compared under a noise level of $\sigma=0.1$, as shown in Figure \ref{fig:FBA_diff_dim}. One can see that while the speed of convergence decreases as the dimensionality increases, the algorithm effectively estimates accurate FOPs as all conditions achieve low training errors within 100 BO iterations. The estimated FOPs generalize well to the testing datasets as the training and testing errors share nearly identical convergence rates and reach the same level after 100 iterations. The ground-truth parameters are recovered in all random instances under different dimensionality, where low parameter errors ($\sim 10^{-4}$) are achieved after 200 iterations.

% spacing control for \midrule

We also record the computation times for solving the 2- and 4-dimensional problems with and without solving the FOPs with 10 parallel workers in a machine on MSI. Table \ref{table:comp_time_parallel} reports the total times as well as the breakdown into the BO and FOPs parts. The BO part includes the time for updating/fitting the GP surrogate and minimizing the acquisition function, whereas the FOPs part accounts for the time required to solve the individual FOPs and compute the loss function value. We observe little differences in the total computation times when the dimension of the problem increases from 2 to 4. Reductions in the total computation time of 46\% and 33\% are achieved in the 2- and 4-dimension problems, respectively, with FOPs solved in parallel. The breakdown times indicate that more computational resources are used to solve the FOPs, and the speed-ups are mainly attributed to the reduced time from the parallel implementation as we see 47$\sim$56\% time reduction in solving the FOPs. Notably, the computation time of the BO part increases around 10$\sim$15\% when the dimension of unknown parameters is doubled.

\vspace{-1em}
\newcommand{\midsepremove}{\aboverulesep = 0mm \belowrulesep = 0mm}
\newcommand{\midsepdefault}{\aboverulesep = 0.605mm \belowrulesep = 0.984mm}
\begin{table}[ht]
\centering

\setlength\extrarowheight{5pt}
\setlength{\tabcolsep}{10pt}
\midsepremove
\begin{threeparttable}
\begin{tabular}{ccccccc}
\toprule[1.5pt]
\multirow{2}{*}{$n_\theta$} & \multirow{2}{*}{Parallelism} & \multicolumn{3}{c}{Median time (s)} & \multicolumn{2}{c}{Reduction (\%)} \\ \cline{3-7} 
                      &                              & Total           & BO           &  FOPs          & Total                    & FOPs                     \\ \midrule[1.5pt]
\multirow{2}{*}{2}                        & No                                               & 3703.0                    & 731.6                    & 2953.7                 & \multirow{2}{*}{45.9}     & \multirow{2}{*}{55.6}    \\ \cline{2-5}
                                          & Yes                                              & 2005.0                    & 698.4                    & 1312.1                 &                           &                          \\ \hline
\multirow{2}{*}{4}                        & No                                               & 3149.4                    & 812.9                    & 2369.5                 & \multirow{2}{*}{33.2}     & \multirow{2}{*}{46.7}    \\ \cline{2-5}
                                          & Yes                                              & 2104.8                    & 825.6                    & 1262.2                 &                           &                          \\ \bottomrule[1.5pt]

\end{tabular}
% \begin{tablenotes}[flushleft]
%     % \footnotesize
%     \item \textit{Notes: Results based on 10 random instances with $|\mathcal{I}|=50$, $ \sigma = 0.1$, and 250 BO iterations. If parallelism is implemented, 50 FOPs are sequentially solved by 10 parallel workers to global optimality. The BO time includes initializing/updating the GP surrogate and minimizing the acquisition function, whereas the FOPs time denotes the time for solving the $|\mathcal{I}|$ FOPs.}
% \end{tablenotes}
\end{threeparttable}
\caption{Computational time of BO4IO as FOPs solved with and without parallel computing. Results based on 10 random instances with $|\mathcal{I}|=50$, $ \sigma = 0.1$, and 250 BO iterations. If parallelism is implemented, 50 FOPs are sequentially solved by 10 parallel workers to global optimality. The BO time includes initializing/updating the GP surrogate and minimizing the acquisition function, whereas the FOPs time denotes the time for solving the $|\mathcal{I}|$ FOPs.}
\label{table:comp_time_parallel}%
\end{table}
\subsubsection{Uncertainty quantification of parameter estimates using profile likelihood}
\begin{figure}[!ht]
\centering
\includegraphics[scale=0.345]{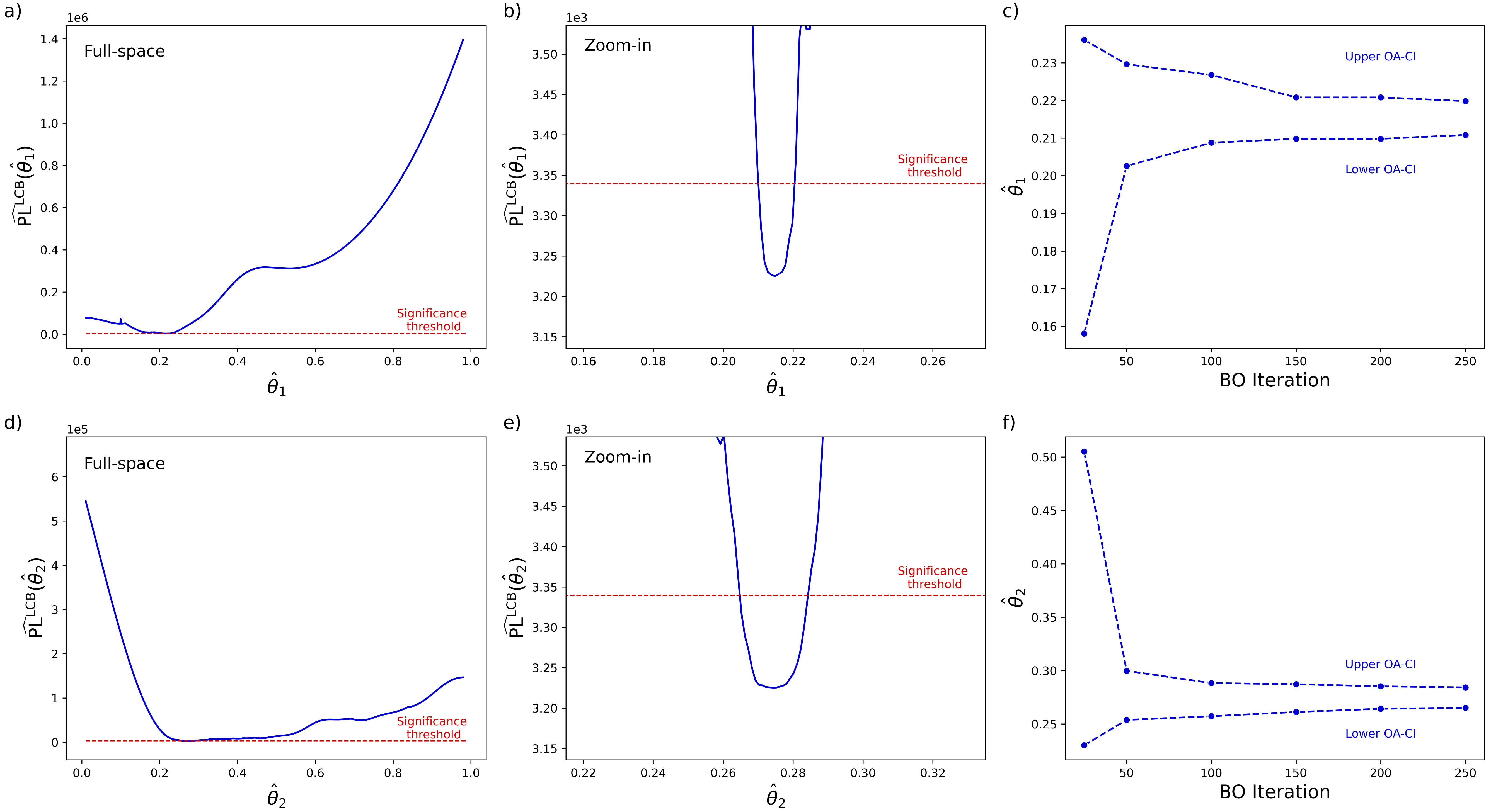}
\caption{Approximation of profile likelihood and identifiability of parameter estimates in a two-dimensional FBA problem ($d=2$). a) and d) show the full-space $\widehat{\text{PL}}^\text{LCB}$ of $\theta_1$ and $\theta_2$, whereas b) and d) are the zoomed-in regions around the significant threshold level as defined in \eqref{eqn:OA-CI}. c) and f) trace the changes in the upper and lower bounds of the OA-CIs of two parameters over the BO iterations. The PL analysis is performed using $\Delta_{\alpha} = \chi^2(0.05, 1)$ and $\rho = 3.84$ (95\% confidence level).}
\label{fig:FBA_PL_2D}
\end{figure}

We next apply the proposed PL analysis to evaluate the uncertainty and identifiability of the parameter estimates of BO4IO in the two-dimensional case. All parameters of the 10 random instances show finite OA-CIs, indicating that they are all structurally identifiable. For simplicity, we choose a random instance from the two-dimensional case as an example. The full-space PL on the two parameters, $\theta_1$ and $\theta_2$, are shown in Figures \ref{fig:FBA_PL_2D}a and \ref{fig:FBA_PL_2D}d, whereas the zoomed-in region are shown in Figures \ref{fig:FBA_PL_2D}b and \ref{fig:FBA_PL_2D}e, respectively. As mentioned, both parameters are identifiable with finite and small OA-CIs. We further evaluate the outer-approximation PL and the approximate CIs across the BO iterations, as shown in Figures \ref{fig:FBA_PL_2D}c and \ref{fig:FBA_PL_2D}f. We observe a monotonic decrease in the width of the OA-CI; the shrinking CI can be attributed to a more accurate GP surrogate of the loss function as the number of BO iterations increases. 
% \subsubsection{Solution quality under varying noise levels}

\begin{figure}[!ht]
\centering
\includegraphics[scale=0.40]{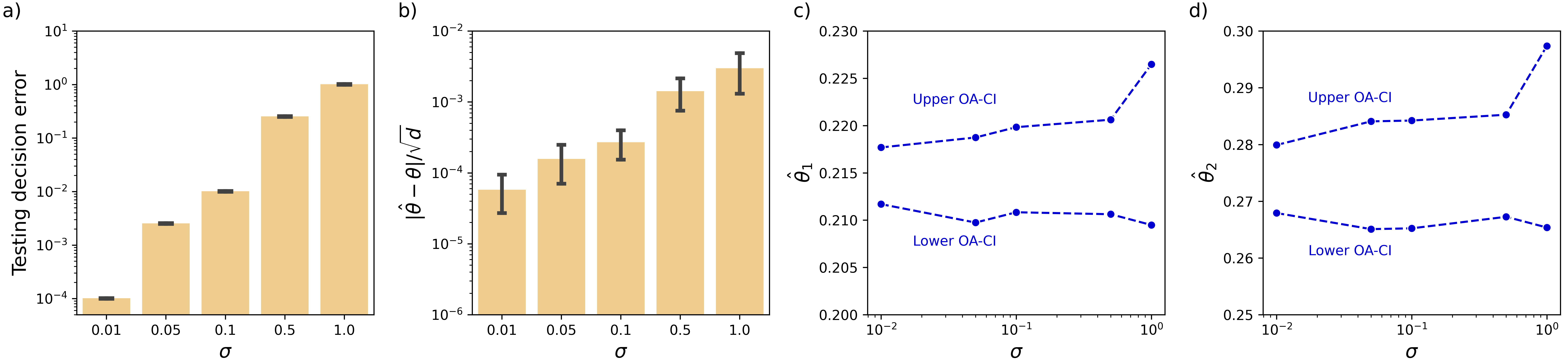}
\caption{Impact of noise levels on the BO4IO performance. (a) Testing decision error and b) parameter error at 250 iterations under varying noise levels. Outer-approximation confidence intervals (OA-CIs) of two parameters, (b) $\theta_1$ and (c) $\theta_2$, under varying noise levels ($\sigma$).  OA-CIs are derived using $\Delta_{\alpha} = \chi^2(0.05, 1)$ and $\rho = 3.84$ (95\% confidence level).}
\label{fig:FBA_PL_sigma}
\end{figure}

Next, again using the two-dimensional case as an example, we test the algorithm performance under varying noise levels ($\sigma$) in the observed decisions (Figure \ref{fig:FBA_PL_sigma}). The average testing decision and $\theta$ error increase as the noise level increases; nonetheless, all conditions converge to the lowest testing errors that can be reached ($\sim \sigma^2$)  at 250 iterations. We further perform the PL analysis to the same random instance discussed in Figure \ref{fig:FBA_PL_2D} at different BO iterations to trace the changes in the parameter CIs over time, as shown in (Figures \ref{fig:FBA_PL_sigma}c and \ref{fig:FBA_PL_sigma}d). While we observe wider OA-CIs under higher noise levels, which aligns with the traditional PL analysis \citep{raue2009structural} where decreasing data quality reduces the parameter identifiability with a wider CI, the estimated parameters remain structurally identifiable under the high noise level.  
\newpage
\subsection{Case study 2: Learning market demands in the standard pooling problem}
In the second case study, we consider a standard pooling problem where an operator blends a set of feedstocks in a pooling network to create different final products that meet desired qualities and demands while minimizing the total cost. Provided below is the formulation of a standard pooling problem \citep{misener2009advances}:
\begin{subequations}\label{eqn:std_pooling}
    \begin{align}
    \minimize_{f,y,z,q} \quad & \sum_{(s,l) \in T_f} \eta_s f_{sl} - \sum_{(l,j)\in T_y} \phi^y_{lj}y_{lj}-\sum_{(s,j)\in T_z} \phi^z_{sj} z_{sj} \\
    \mathrm {subject \; to} \quad &  \sum_{l:(s,l) \in T_f} f_{sl} + \sum_{j:(s,j)\in T_z}z_{sj} \leq A^U_s \quad \forall s \in \mathcal{S} \label{eqn:std_pooling-b}\\
    \quad &  \sum_{s:(s,l) \in T_f} f_{sl} - \sum_{j:(l,j) \in T_y}y_{lj} = 0 \quad \forall l \in \mathcal{L} \label{eqn:std_pooling-c}\\
    \quad &  \sum_{s:(s,l) \in T_f} C_{sk} f_{sl} - p_{lk} \sum_{j:(l,j) \in T_y} y_{lj} =0 \quad \forall l \in \mathcal{L},\ k \in \mathcal{K} \label{eqn:std_pooling-d}\\
    \quad & \sum_{l:(l,j) \in T_y} p_{lk} y_{lj} + \sum_{s:(s,j) \in T_z}C_s z_{sj} \leq P^U_{jk} \left(\sum_{l:(l,j) \in T_y} y_{lj} + \sum_{s:(s,j) \in T_z} z_{sj} \right) \quad \forall j \in \mathcal{J},\ k \in \mathcal{K} \label{eqn:std_pooling-e}\\
    \quad & \ \sum_{l:(l,j) \in T_y} y_{lj} + \sum_{s:(s,j) \in T_z} z_{sj} \leq \theta_j \quad \forall j \in \mathcal{J} \label{eqn:std_pooling-f}\\
    \quad & f_{sl} \geq 0 \quad \forall s \in \mathcal{S},\ l \in \mathcal{L}\label{eqn:std_pooling-g}\\
    \quad & y_{lj} \geq 0 \quad \forall l \in \mathcal{L},\ j \in \mathcal{J} \label{eqn:std_pooling-h}\\
    \quad & z_{sj} \geq 0 \quad \forall s \in \mathcal{S},\ j \in \mathcal{J} \label{eqn:std_pooling-i}\\
    \quad & p_{lk} \geq 0 \quad \forall l \in \mathcal{L},\ k \in \mathcal{K}, \label{eqn:std_pooling-j}
    \end{align}
\end{subequations}
where $\mathcal{S}$, $\mathcal{L}$, $\mathcal{J}$, and $\mathcal{K}$ are the sets of input feedstocks, mixing pools, output products, and quality attributes, respectively. As incoming feedstocks can connect to a pool or directly to an output, sets $T_f$, $T_y$, and $T_z$ denote the existing streams from input $s$ to pool $l$, pool $l$ to output $j$, and input $s$ to output $j$, respectively. The quality per unit of feedstock $s$ of attribute $k$ is denoted as $C_{sk}$. The cost per unit of feedstock $s$ is denoted as $\eta_s$. The revenue per unit flow from pool $l$ to output $j$ and input $s$ to output $j$ are denoted by $\phi^y_{ij}$ and $\phi^z_{sj}$, respectively. The upper limit of the quality attribute $k$ in each output product $j$ is noted as $P^U_{jk}$. Decision variables include $f_{sl}$, $y_{lj}$, and $z_{sj}$ denoting the flow from input $s$ to pool $l$, pool $l$ to product $j$, and input $s$ to output $j$, respectively, whereas the quality level in pool $l$ of attribute $k$ is denoted by $p_{lk}$.

In problem \eqref{eqn:std_pooling}, we assume that each feedstock $s$ has limited availability $A^U_s$ as indicated in constraints \eqref{eqn:std_pooling-b}; the material and quality balance at pool $l$ are maintained through constraints \eqref{eqn:std_pooling-c} and \eqref{eqn:std_pooling-d}, respectively; the upper acceptable product quality constraint is set in \eqref{eqn:std_pooling-e}. We consider a scenario in which the demand for each product $j$ ($\theta_j$ in constraints \eqref{eqn:std_pooling-f}) is unknown as the unknown dimension denoted as $d$. Notably, constraints \eqref{eqn:std_pooling-e} contain bilinear terms that render the overall optimization problem nonconvex. The goal is to apply BO4IO to estimate the $\theta_j$ values by observing a part of the decisions, namely $f_i$ and $y_i$, based on varying input conditions $u_i=(A^U, P^U, \eta, \phi^y,\phi^z)_i$ in a set of observations $\mathcal{I}$. 

We test the proposed BO4IO framework on two benchmark pooling networks, Haverly1 \citep{Haverly1978studies} and Foulds3 \citep{foulds1992bilinear} where their network specifications ($|\mathcal{S}|$,$|\mathcal{L}|$,$|\mathcal{J}|$,$|\mathcal{K}|$) are (3,1,2,1) and (32,8,16,1), respectively. For each random instance of the IOP, the synthetic dataset is generated as follows. We first generate a set of ground-truth $\theta_j\sim \mathcal{U}(0.5,1.0)$ for every $j\in \mathcal{J}$ with randomized feedstock price $\eta\sim\mathcal{U}(\bar{\eta}^\text{min},\bar{\eta}^\text{max})$ and product revenue $\phi{lj}^y, \phi{sj}^z \sim \mathcal{U}(\bar{\phi}^\text{min},\bar{\phi}^\text{max})$, where $\bar{\eta}^{\text{min}/\text{max}}$ and $\bar{\phi}^{\text{min}/\text{max}}$ denote the minimum or maximum of the nominal $\eta$ and $\phi$ values in the original problems. Using the same values of ($\theta, \eta, \phi^y, \phi^z$) and the nominal values of $P^U$, we solved \eqref{eqn:std_pooling} with randomized availability $A^U \sim \mathcal{U}(0.5,1.0)$ to collect the true optimal solutions $x_i^*=(f_i^*,y_i^*)$ for every $i \in \mathcal{I}$. Lastly, we generated noisy observations of decisions $x_i=x_i^*+\gamma$, where $\gamma \sim \mathcal{N}(0,\sigma^2\mathbb{I})$ for all $i \in \mathcal{I}$. A separate set of testing data is generated based on the same randomization procedure with 50 experiments. The computational statistics are obtained from the results of 10 random instances of each IOP.

\subsubsection{Robust performance against network sizes}
\begin{figure}[!hb]
\centering
\includegraphics[scale=0.55]{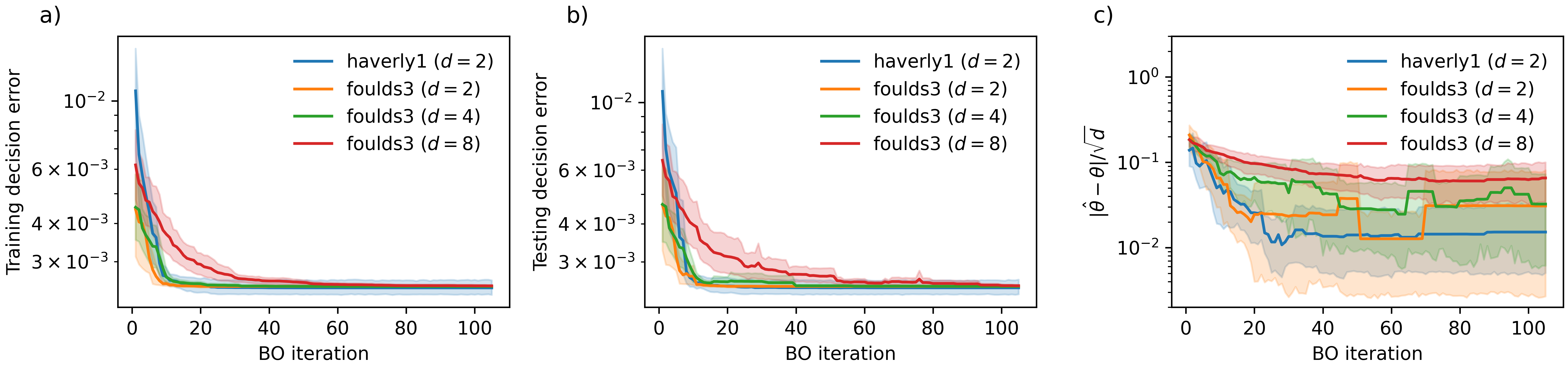}
\caption{Effect of network size to the BO4IO performance. Convergence analysis of (a) training error, (b) testing error, and (c) parameter error with two benchmark problems, as the numbers shown in the parentheses in the legend denote the varying dimensionality $d$. Training and testing loss are defined as the average standardized decision error, $\sum_{i\in\mathcal{I}}(x_i-\hat{x}_i)^\top (x_i-\hat{x}_i)/|\mathcal{I}|/|\mathcal{R}|$ and calculated based on the training and testing datasets, respectively, whereas parameter error denotes the difference between the ground-truth ($\theta$) and estimated ($\hat{\theta}^*$) values. Here, the solid lines and shaded areas respectively denote the medians and confidence intervals of the corresponding loss across the 10 random instances. The synthetic dataset of each random instance is generated under the setting of $|\mathcal{I}| = 50$ and $\sigma=0.05$.}
\label{fig:std_diff_prob}
\end{figure}

We first test how robust the algorithm performance is under different pooling network sizes by applying it to the two benchmark problems with $|\mathcal{I}|=50$ and $\sigma=0.05$. To provide a fair comparison between the two problems, we apply the algorithm to learn different dimensionality $d$ of the unknown demands $\theta$ in the Foulds3 problem. We trace the progressions of training, testing, and parameter errors over the BO iterations, as shown in Figure \ref{fig:std_diff_prob}. The algorithm efficiently learns the unknown parameters regardless of the network complexity, as we observe a similar convergence speed in training and testing errors in the two network problems under the same number of unknown parameters. Consistent with the first case study, a lower convergence rate is observed when the dimensionality increases in the Foulds3 problem. It should be noted that Foulds3 shows a higher parameter error than Haverly1, and the error increases with the dimensionality. We further perform the proposed PL analysis and find that this is due to the larger numbers of non-identifiable parameters in the Foulds3 problems (data not shown). 
\vspace{-0.05em}

\subsubsection{Data efficiency and identifiability analysis}
\begin{figure}[!ht]
\centering
\includegraphics[scale=0.2]{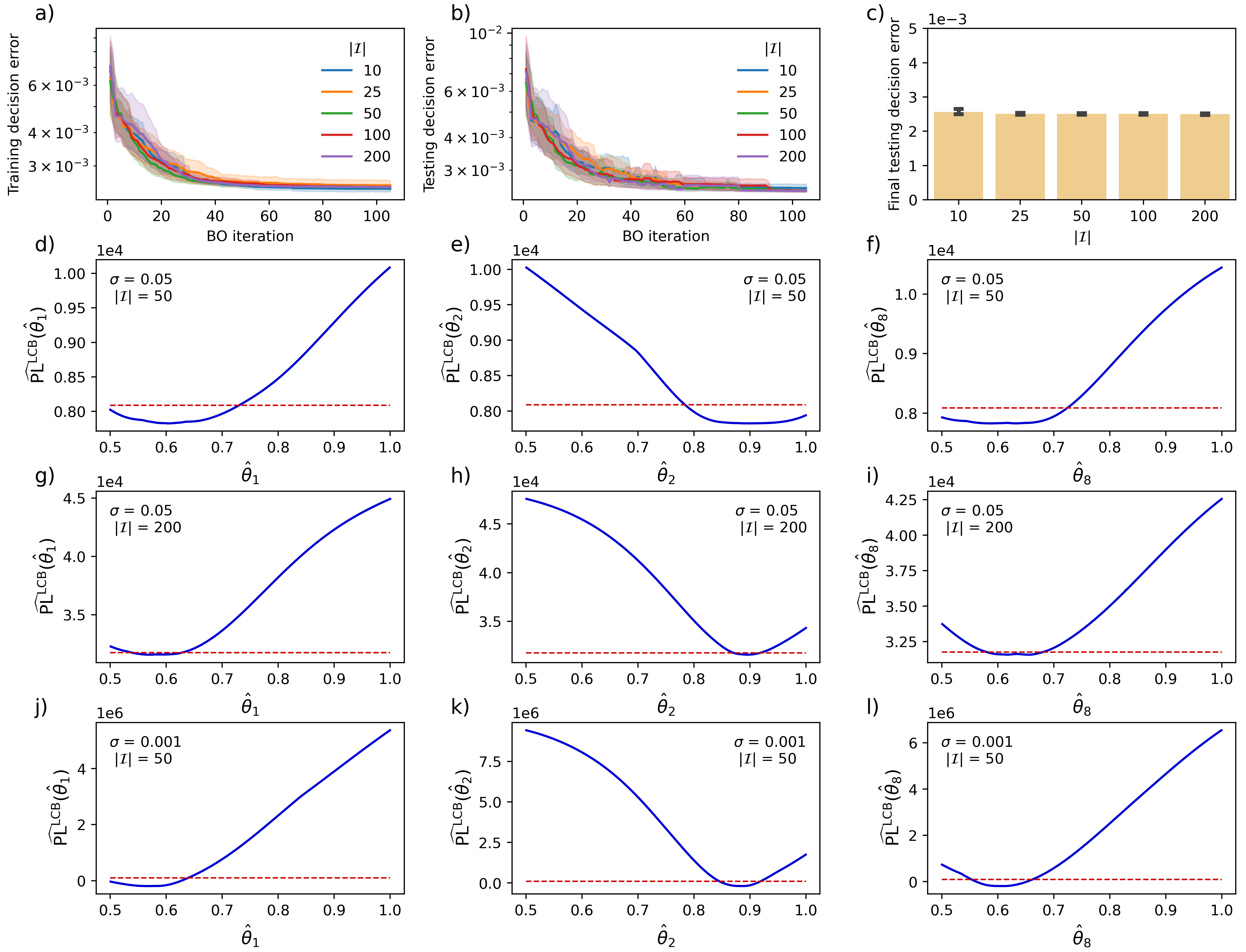}
\caption{Effect of training dataset size on the BO4IO performance using Foulds3 with $d=8$ and $\sigma=0.05$ as an example. Convergence analysis of (a) training error and (b) testing error as well as (c) final testing error at 100 iterations under different sizes of training datasets. d)-f) show the approximate PL of three practically non-identifiable parameters in the selected random IOP instance under $|\mathcal{I}| = 50$ and $\sigma=0.05$. g)-i) show the approximate PL of the same three parameters as d)-f) estimated under $|\mathcal{I}| = 200$ and $\sigma=0.05$. j)-l) show the approximate PL of the same three parameters as d)-f) estimated under $|\mathcal{I}| = 50$ and $\sigma=0.001$. In a) and b), the solid lines and shaded areas respectively denote the medians and confidence intervals of the errors across the 10 random instances. In d)-l), the red dashed line denotes the significance threshold for the identifiability analysis. The PL analysis is performed using $\Delta_{\alpha} = \chi^2(0.05, 1)$ and $\rho = 3.84$ (95\% confidence level).}
\label{fig:std_diff_data_size}
\end{figure}

We next test the data efficiency of the algorithm by learning 8 unknown demands in Foulds3 based on different sizes of the training datasets, namely $|\mathcal{I}|=\{10,25,50,100,200\}$. The convergence profiles for the training and testing errors under different conditions are shown in Figures \ref{fig:std_diff_data_size}a and \ref{fig:std_diff_data_size}b. One can see that BO4IO requires only few data points to achieve a low decision error as there is little difference in the convergence rates for the different training dataset sizes and no significant difference in the final testing error as shown in Figure \ref{fig:std_diff_data_size}c.

We further perform the PL analysis to evaluate the impact of training dataset size on parameter identifiability. For simplicity, we select one random IOP instance as an example. The instance shows 5 structurally identifiable parameters (data not shown) and 3 practically non-identifiable parameters (Figure \ref{fig:std_diff_data_size}d-f) when learned from a training dataset with $|\mathcal{I}|=50$ and $\sigma=0.05$. Practical non-identifiability is typically due to insufficient data quality or quantity \citep{raue2009structural}. One can see that most of the three practically non-identifiable parameters become structurally identifiable when we further increase the number of data points to $|\mathcal{I}|=200$ (Figure \ref{fig:std_diff_data_size}g-i) or reduce the observation noise to $\sigma=0.001$ ($|\mathcal{I}|=200$ Figure \ref{fig:std_diff_data_size}j-l). Notably, the OA-CIs of the 5 other parameters maintain their structural identifiability with smaller OA-CIs (data not shown). The proposed PL analysis provides a means of uncertainty quantification on the BO4IO parameter estimates, and this example demonstrates its importance in assessing the identifiability of the unknown parameters.

\subsection{Case study 3: Learning quality requirements in the generalized pooling problem}
In the third case study, we apply the BO4IO algorithm to learn the product quality parameters in the constraints of a generalized pooling problem \citep{misener2009advances}, which is an extension of \eqref{eqn:std_pooling} with additional discrete variables. The problem can be formulated as the following nonconvex mixed-integer nonlinear program: 
\begin{subequations}\label{eqn:gen_pool}
    \begin{align}
    \minimize_{f,y,z,q,\gamma^\text{pool},\gamma^\text{init}} \quad & \sum_{(s,l) \in T_f} (\eta_s + \eta^{f}_{sl})f_{sl} - \sum_{(l,j)\in T_y} \phi^y_{lj}y_{lj}-\sum_{(s,j)\in T_z} \phi^z_{sj} z_{sj} + \sum_{l\in\mathcal{L}}\eta^{\text{pool}}_l\gamma^\text{pool}_l + \sum_{s\in\mathcal{S}}\eta^{\text{init}}_s\gamma^\text{init}_s\\
    \mathrm {subject \; to} \quad &  \sum_{l:(s,l) \in T_f} f_{sl} + \sum_{j:(s,j)\in T_z}z_{sj} \leq A^U_s \gamma^\text{init}_s \quad \forall s \in \mathcal{S}\\
    \quad &  \sum_{s:(s,l) \in T_f} f_{sl} - \sum_{j:(l,j) \in T_y}y_{lj} = 0 \quad \forall l \in \mathcal{L} \\
    \quad &  \sum_{s:(s,l) \in T_f} C_{sk} f_{sl} - p_{lk} \sum_{j:(l,j) \in T_y} y_{lj} =0 \quad \forall l \in \mathcal{L},\ k \in \mathcal{K} \\
    \quad & \sum_{l:(l,j) \in T_y} p_{lk} y_{lj} + \sum_{s:(s,j) \in T_z}C_{sk} z_{sj} \leq \theta_{jk} (\sum_{l:(l,j) \in T_y} y_{lj} + \sum_{s:(s,j) \in T_z} z_{sj}) \quad \forall j \in \mathcal{J},\ k \in \mathcal{K}  \label{eqn:gen_pool-e} \\
    \quad & \sum_{l:(l,j) \in T_y} y_{lj} + \sum_{s:(s,j) \in T_z} z_{sj} =D_j \quad \forall j \in \mathcal{J} \label{eqn:gen_pool-f} \\
    \quad & \sum_{j:(l,j)\in T_y} y_{lj} \leq S_l \gamma^\text{pool}_l \quad \forall l \in \mathcal{L} \label{eqn:gen_pool-g} \\
    \quad & \sum_{l:(s,l)\in T_f} f_{sl} \leq A_s^U \gamma^\text{init}_s \quad \forall s \in \mathcal{S} \label{eqn:gen_pool-h} \\
    \quad & f_{sl} \geq 0 \quad \forall s \in \mathcal{S},\ l \in \mathcal{L}\\
    \quad & y_{lj} \geq 0 \quad \forall l \in \mathcal{L},\ j \in \mathcal{J}\\
    \quad & z_{sj} \geq 0 \quad \forall s \in \mathcal{S},\ j \in \mathcal{J}\\
    \quad & p_{lk} \geq 0 \quad \forall l \in \mathcal{L},\ k \in \mathcal{K}\\
    \quad & \gamma^{\text{init}}_s \in \{0,1\} \quad \forall s \in \mathcal{S} \\
    \quad & \gamma^{\text{pool}}_l \in \{0,1\} \quad \forall l \in \mathcal{L}.
    \end{align}
\end{subequations}
The model is an extension to \eqref{eqn:std_pooling} with some minor changes. Binary variables are introduced to consider the installation decisions of feedstock $s$ and pool $l$, denoted as $\gamma^{\text{init}}_s$ and $\gamma^{\text{pool}}_l$, respectively. In addition, new parameters are added to consider the installation cost of feedstock $s$ and pool $l$, denoted as $\eta^\text{init}_s$ and $\eta^\text{pool}_l$, whereas $\eta^f_{sl}$ is the cost per unit flow from feedstock $s$ to pool $l$. Two additional constraints \eqref{eqn:gen_pool-g} and \eqref{eqn:gen_pool-h} are included to constrain the flow variables accounting for the binary feedstock and pool installation, where $S_l$ denotes the pool capacity of pool $l$. 

We consider a benchmark problem called Lee1 \citep{lee2003global}, which has a network specification ($|\mathcal{S}|$,$|\mathcal{L}|$,$|\mathcal{J}|$,$|\mathcal{K}|$) of (5,4,3,2).
This study addresses the scenario in which each product stream $j$ has a demand that needs to be satisfied ($D_j$ in \eqref{eqn:gen_pool-f}) with unknown quality requirements ($\theta_{jk}$ in \eqref{eqn:gen_pool-e}). We assume the the quality requirements $\theta_{jk}$ of all attributes sum up to one such that $\sum_{k\in\mathcal{K}}\theta_{j,k}=1$ for every $j\in \mathcal{J}$. Since Lee1 contains two quality attributes, we only need to estimate one of the quality attributes, and the number of unknown parameters thus reduces to $|\mathcal{J}| = 3$. The synthetic data are generated based on the following procedure. For each instance of IOP, we first generate the ground-truth quality requirements $\theta_{j1} \sim \mathcal{U}(0.2,0.6)$ and randomized product revenue $\phi_{lj}^y, \phi_{sj}^z \sim \mathcal{U}(0.5\bar{\phi}^\text{min},1.5\bar{\phi}^\text{max})$, where $\bar{\phi}^{\text{min}/\text{max}}$ denote the minimum or maximum of the nominal $\phi$ values in the original problem. Using the fixed $\theta$, $\phi^y$, $\phi^z$ values, we then solve the $|\mathcal{I}|$ FOPs with randomized $A^U_s, D_j \sim \mathcal{U}(0.5,1.5)$. Similar to the second case study, only a subset of decisions is collected to generate the noisy observations, where $x_i^* = f_i^*$ and $x_i=x_i^*+\gamma$ with $\gamma \sim \mathcal{N}(0,\sigma^2\mathbb{I})$ for all $i \in \mathcal{I}$. A separate testing dataset of 50 experiments is generated to validate the model prediction.

\begin{figure}[!ht]
\centering
\includegraphics[scale=0.54]{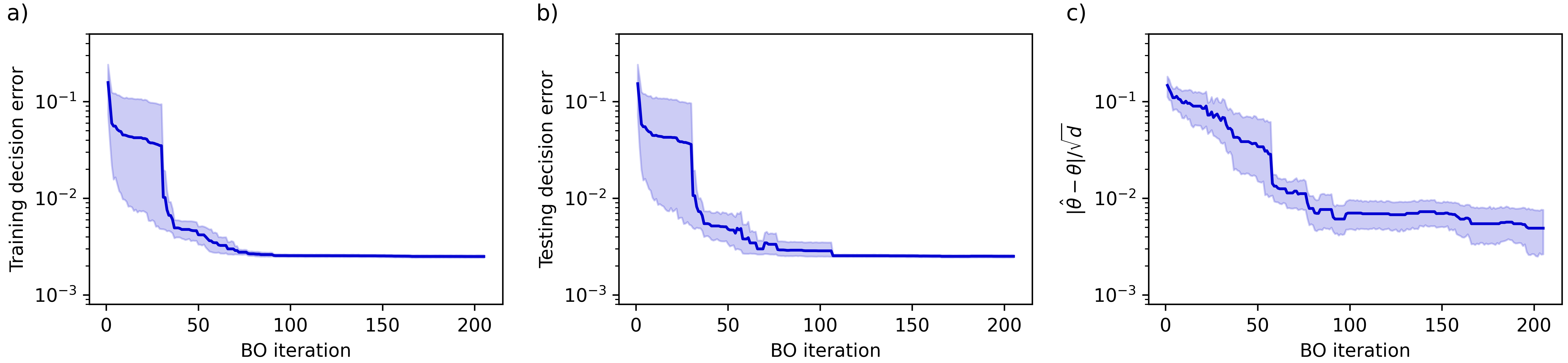}
\caption{Algorithm performance for the generalized pooling problem Lee1. Convergence analysis of (a) training error, (b) testing error, and (c) parameter error with three unknown parameters. Training and testing errors are defined as the average standardized decision errors, $\sum_{i\in\mathcal{I}}(x_i-\hat{x}_i)^\top (x_i-\hat{x}_i)/|\mathcal{I}|/|\mathcal{R}|$ and calculated based on the training and testing datasets, respectively, whereas parameter error denotes the difference between the ground-truth ($\theta$) and estimated ($\hat{\theta}^*$) values. Here, the solid lines and shaded areas respectively denote the medians and confidence intervals of the corresponding loss across the 10 random instances. The synthetic dataset of each random instance is generated under the setting of $|\mathcal{I}| = 50$ and $\sigma=0.05$.}
\label{fig:gen_conv}
\end{figure}

\begin{figure}[!ht]
\centering
\includegraphics[scale=0.55]{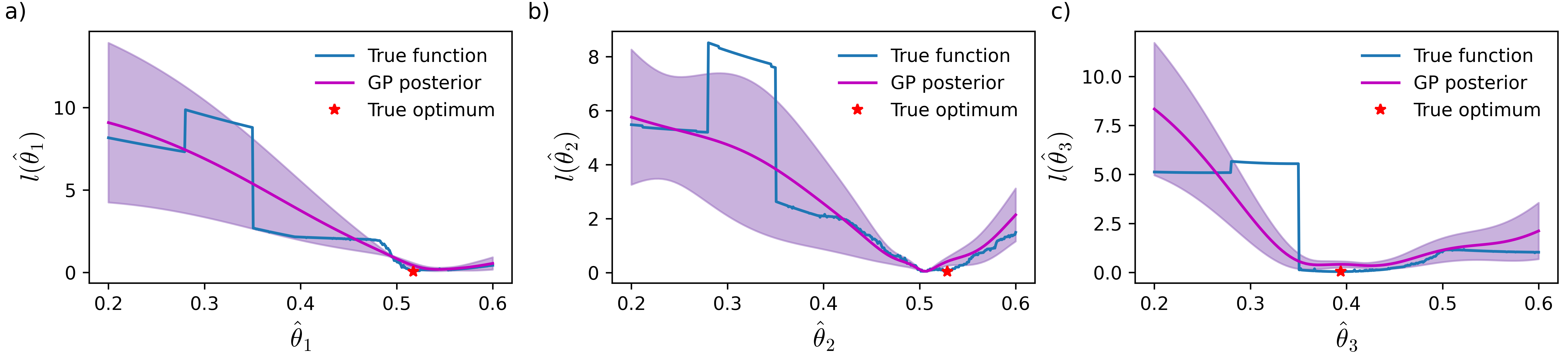}
\caption{Comparison of the changes in the true loss function value (blue curve) and the predicted GP posterior around the best-found solution in the generalized pooling problem (Lee1). The purple lines and shaded areas denote the means and confidence intervals of the predicted GP posteriors.}
\label{fig:gen_SI}
\end{figure}

\begin{figure}[!ht]
\centering
\includegraphics[scale=0.6]{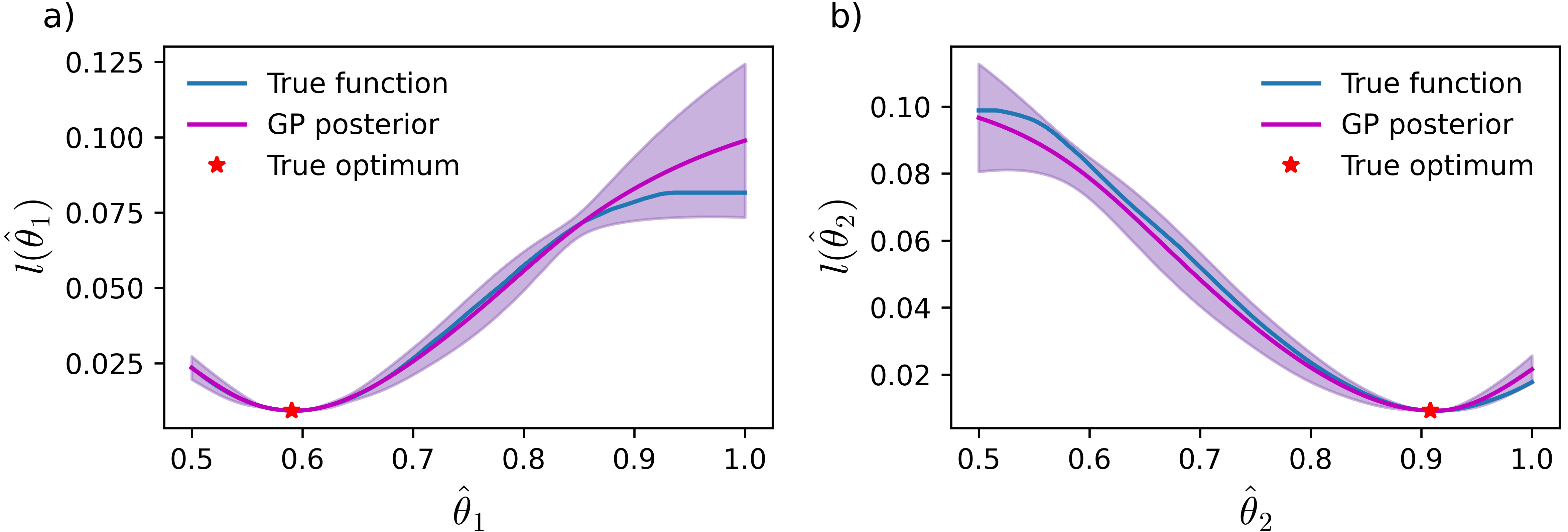}
\caption{Comparison of the changes in the true loss function value (blue curve) and the predicted GP posterior around the best-found solution in the standard pooling problem (Haverly1). The purple lines and shaded areas denote the means and confidence intervals of the predicted GP posteriors.}
\label{fig:gen_SI_Std}
\end{figure}

We apply the algorithm to 10 random instances with $|\mathcal{I}|=50$ and $\sigma = 0.05$. The convergence plots are shown in Figure \ref{fig:gen_conv}. The training and testing decision errors converge to a low level within 200 iterations. Compared to the standard pooling problem with a similar specification and number of unknowns (e.g. Haverly1), the algorithm converges more slowly and achieves a higher parameter error in the generalized pooling problem. To further investigate the cause of the slow convergence, we plot the true loss function as a function of each unknown parameter (while keeping the other parameters at their best-found values) for Lee1 and Haverly1, and compare it with the predicted GP posterior, as shown in Figures \ref{fig:gen_SI} and \ref{fig:gen_SI_Std}. In the generalized pooling case (Figure \ref{fig:gen_SI}), we observe discrete (step) changes in the loss function values, and the confidence intervals of the GP posterior do not fully cover the true function. Despite the nonsmoothness of the true loss function, the GP model provides a continuous approximation and leads to a parameter estimate that is close to the true minimum, albeit with a slightly higher parameter error compared to the standard pooling problem. In contrast, the GP posterior provides a very good approximation to the true loss function in the standard pooling problem (Figure \ref{fig:gen_SI_Std}) since the loss function is smooth, and the algorithm finds a solution very close to the global optimum. This analysis indicates that the slow convergence in the generalized pooling problem could be attributed to the discrete nature of the mixed-integer problem, where the poor GP fitness to the discontinuous and nonsmooth region can potentially hinder the algorithm. One may consider nonstationary kernels \citep{paciorek2003nonstationary} or other surrogate types, e.g. deep GP or Bayesian neural networks \citep{damianou2013deep,springenberg2016bayesian}, to capture the nonsmoothness and discontinuity of the target function, which we plan to explore more in future work. 

\section{Conclusions} \label{sec:conclusions}
In this work, we developed a Bayesian optimization approach for data-driven inverse optimization problem, which aims to infer unknown parameters of an optimization model from observed decisions and is commonly formulated as a large-scale bilevel optimization problem. In the proposed  BO4IO framework, the IO loss function is treated as a black box and approximated with a Gaussian process model. At each iteration of the algorithm, the loss function is evaluated and the next candidate solution is determined by minimizing an acquisition function based on the GP posterior. The key advantage of BO4IO is that the evaluation of the loss function can be achieved by directly solving the forward optimization problem with the current parameter estimates for every data point, which circumvents the need for a complex single-level reformulation of the bilevel program or a special, e.g. cutting-plane-based, algorithm. Consequently, BO4IO remains applicable even when the FOP is nonlinear and nonconvex, and it can consider both continuous and discrete decision variables. Furthermore, we leverage the GP posterior to derive confidence bounds on the profile likelihood to compute inner and outer approximations of the confidence interval on each parameter estimate, which provides a means of uncertainty quantification and assessing the identifiability of the unknown parameters.

We tested the proposed algorithm and PL method in three computational case studies. The extensive computational results demonstrate the efficacy of BO4IO in estimating unknown parameters in various classes of FOPs, ranging from convex nonlinear to nonconvex mixed-integer nonlinear programs. In the tested instances, the algorithm proved to be efficient in finding parameter estimates with very low prediction errors in relatively few ($\sim$100) BO iterations. Using the proposed PL analysis, we were able to determine the different extents of identifiability for the unknown parameters and also demonstrated how the identifiability of a parameter can change with the quantity and quality of data.

Finally, we would like to highlight that several important directions for future work remain, including integrating high-dimensional BO techniques \citep{snoek2015scalable, eriksson2021high} for learning large sets of parameter and developing a customized surrogate for mixed-integer problems to further improve the algorithm performance.

\section*{Acknowledgement}
The authors thank Vikram Kumar, who was an undergraduate student at the University of Minnesota at the time, for conducting preliminary computational results for this work. The authors gratefully acknowledge the support from the National Science Foundation under Grants \#2044077 and \#2237616. Computational resources were provided by the Minnesota Supercomputing Institute at the University of Minnesota.

\singlespacing

\bibliographystyle{plainnat}
\bibliography{References.bib}

\begin{thebibliography}{90}
\providecommand{\natexlab}[1]{#1}
\providecommand{\url}[1]{\texttt{#1}}
\expandafter\ifx\csname urlstyle\endcsname\relax
  \providecommand{\doi}[1]{doi: #1}\else
  \providecommand{\doi}{doi: \begingroup \urlstyle{rm}\Url}\fi

\bibitem[Ahuja and Orlin(2001)]{ahuja2001inverse}
Ravindra~K Ahuja and James~B Orlin.
\newblock Inverse optimization.
\newblock \emph{Operations research}, 49\penalty0 (5):\penalty0 771--783, 2001.

\bibitem[Ajayi et~al.(2022)Ajayi, Lee, and Schaefer]{ajayi2022objective}
Temitayo Ajayi, Taewoo Lee, and Andrew~J Schaefer.
\newblock Objective selection for cancer treatment: An inverse optimization
  approach.
\newblock \emph{Operations Research}, 70\penalty0 (3):\penalty0 1717--1738,
  2022.

\bibitem[Aswani et~al.(2018)Aswani, Shen, and Siddiq]{aswani2018inverse}
Anil Aswani, Zuo-Jun Shen, and Auyon Siddiq.
\newblock Inverse optimization with noisy data.
\newblock \emph{Operations Research}, 66\penalty0 (3):\penalty0 870--892, 2018.

\bibitem[Auer(2002)]{auer2002using}
Peter Auer.
\newblock Using confidence bounds for exploitation-exploration trade-offs.
\newblock \emph{Journal of Machine Learning Research}, 3\penalty0
  (Nov):\penalty0 397--422, 2002.

\bibitem[Balandat et~al.(2020)Balandat, Karrer, Jiang, Daulton, Letham, Wilson,
  and Bakshy]{balandat2020botorch}
Maximilian Balandat, Brian Karrer, Daniel Jiang, Samuel Daulton, Ben Letham,
  Andrew~G Wilson, and Eytan Bakshy.
\newblock Botorch: A framework for efficient monte-carlo bayesian optimization.
\newblock \emph{Advances in neural information processing systems},
  33:\penalty0 21524--21538, 2020.

\bibitem[Balsa-Canto et~al.(2008)Balsa-Canto, Alonso, and
  Banga]{balsa2008computational}
Eva Balsa-Canto, Antonio~A Alonso, and Julio~R Banga.
\newblock Computational procedures for optimal experimental design in
  biological systems.
\newblock \emph{IET systems biology}, 2\penalty0 (4):\penalty0 163--172, 2008.

\bibitem[Bandara et~al.(2009)Bandara, Schl{\"o}der, Eils, Bock, and
  Meyer]{bandara2009optimal}
Samuel Bandara, Johannes~P Schl{\"o}der, Roland Eils, Hans~Georg Bock, and
  Tobias Meyer.
\newblock Optimal experimental design for parameter estimation of a cell
  signaling model.
\newblock \emph{PLoS computational biology}, 5\penalty0 (11):\penalty0
  e1000558, 2009.

\bibitem[Bellman and {\AA}str{\"o}m(1970)]{bellman1970structural}
Ror Bellman and Karl~Johan {\AA}str{\"o}m.
\newblock On structural identifiability.
\newblock \emph{Mathematical biosciences}, 7\penalty0 (3-4):\penalty0 329--339,
  1970.

\bibitem[Bergstra et~al.(2011)Bergstra, Bardenet, Bengio, and
  K{\'e}gl]{bergstra2011algorithms}
James Bergstra, R{\'e}mi Bardenet, Yoshua Bengio, and Bal{\'a}zs K{\'e}gl.
\newblock Algorithms for hyper-parameter optimization.
\newblock \emph{Advances in neural information processing systems}, 24, 2011.

\bibitem[Berkenkamp et~al.(2019)Berkenkamp, Schoellig, and
  Krause]{berkenkamp2019no}
Felix Berkenkamp, Angela~P Schoellig, and Andreas Krause.
\newblock No-regret bayesian optimization with unknown hyperparameters.
\newblock \emph{Journal of Machine Learning Research}, 20\penalty0
  (50):\penalty0 1--24, 2019.

\bibitem[Bertsimas et~al.(2015)Bertsimas, Gupta, and
  Paschalidis]{bertsimas2015data}
Dimitris Bertsimas, Vishal Gupta, and Ioannis~Ch Paschalidis.
\newblock Data-driven estimation in equilibrium using inverse optimization.
\newblock \emph{Mathematical Programming}, 153:\penalty0 595--633, 2015.

\bibitem[Beykal et~al.(2020)Beykal, Avraamidou, Pistikopoulos, Onel, and
  Pistikopoulos]{beykal2020domino}
Burcu Beykal, Styliani Avraamidou, Ioannis~PE Pistikopoulos, Melis Onel, and
  Efstratios~N Pistikopoulos.
\newblock Domino: Data-driven optimization of bi-level mixed-integer nonlinear
  problems.
\newblock \emph{Journal of Global Optimization}, 78:\penalty0 1--36, 2020.

\bibitem[Bonarius et~al.(1996)Bonarius, Hatzimanikatis, Meesters, De~Gooijer,
  Schmid, and Tramper]{bonarius1996metabolic}
Hendrik~PJ Bonarius, Vassily Hatzimanikatis, Koen~PH Meesters, Cornelis~D
  De~Gooijer, Georg Schmid, and Johannes Tramper.
\newblock Metabolic flux analysis of hybridoma cells in different culture media
  using mass balances.
\newblock \emph{Biotechnology and bioengineering}, 50\penalty0 (3):\penalty0
  299--318, 1996.

\bibitem[Burgard and Maranas(2003)]{burgard2003optimization}
Anthony~P Burgard and Costas~D Maranas.
\newblock Optimization-based framework for inferring and testing hypothesized
  metabolic objective functions.
\newblock \emph{Biotechnology and bioengineering}, 82\penalty0 (6):\penalty0
  670--677, 2003.

\bibitem[Burton and Toint(1992)]{burton1992instance}
Didier Burton and Ph~L Toint.
\newblock On an instance of the inverse shortest paths problem.
\newblock \emph{Mathematical programming}, 53:\penalty0 45--61, 1992.

\bibitem[Bynum et~al.(2021)Bynum, Hackebeil, Hart, Laird, Nicholson, Siirola,
  Watson, and Woodruff]{bynum2021pyomo}
Michael~L. Bynum, Gabriel~A. Hackebeil, William~E. Hart, Carl~D. Laird,
  Bethany~L. Nicholson, John~D. Siirola, Jean-Paul Watson, and David~L.
  Woodruff.
\newblock \emph{Pyomo--optimization modeling in python}, volume~67.
\newblock Springer Science \& Business Media, third edition, 2021.

\bibitem[Chan and Kaw(2020)]{chan2020inverse}
Timothy~CY Chan and Neal Kaw.
\newblock Inverse optimization for the recovery of constraint parameters.
\newblock \emph{European Journal of Operational Research}, 282\penalty0
  (2):\penalty0 415--427, 2020.

\bibitem[Chan et~al.(2014)Chan, Craig, Lee, and Sharpe]{chan2014generalized}
Timothy~CY Chan, Tim Craig, Taewoo Lee, and Michael~B Sharpe.
\newblock Generalized inverse multiobjective optimization with application to
  cancer therapy.
\newblock \emph{Operations Research}, 62\penalty0 (3):\penalty0 680--695, 2014.

\bibitem[Chan et~al.(2023)Chan, Mahmood, and Zhu]{chan2023inverse}
Timothy~CY Chan, Rafid Mahmood, and Ian~Yihang Zhu.
\newblock Inverse optimization: Theory and applications.
\newblock \emph{Operations Research}, 2023.

\bibitem[Choon et~al.(2014)Choon, Mohamad, Deris, Illias, Chong, Chai, Omatu,
  and Corchado]{choon2014differential}
Yee~Wen Choon, Mohd~Saberi Mohamad, Safaai Deris, Rosli~Md Illias, Chuii~Khim
  Chong, Lian~En Chai, Sigeru Omatu, and Juan~Manuel Corchado.
\newblock Differential bees flux balance analysis with optknock for in silico
  microbial strains optimization.
\newblock \emph{PloS one}, 9\penalty0 (7):\penalty0 e102744, 2014.

\bibitem[Cobelli and Distefano~III(1980)]{cobelli1980parameter}
Claudio Cobelli and Joseph~J Distefano~III.
\newblock Parameter and structural identifiability concepts and ambiguities: a
  critical review and analysis.
\newblock \emph{American Journal of Physiology-Regulatory, Integrative and
  Comparative Physiology}, 239\penalty0 (1):\penalty0 R7--R24, 1980.

\bibitem[Damianou and Lawrence(2013)]{damianou2013deep}
Andreas Damianou and Neil~D Lawrence.
\newblock Deep gaussian processes.
\newblock In \emph{Artificial intelligence and statistics}, pages 207--215.
  PMLR, 2013.

\bibitem[Dogan and Prestwich(2023)]{dogan2023bilevel}
Vedat Dogan and Steven Prestwich.
\newblock Bilevel optimization by conditional bayesian optimization.
\newblock In \emph{International Conference on Machine Learning, Optimization,
  and Data Science}, pages 243--258. Springer, 2023.

\bibitem[Eriksson and Jankowiak(2021)]{eriksson2021high}
David Eriksson and Martin Jankowiak.
\newblock High-dimensional bayesian optimization with sparse axis-aligned
  subspaces.
\newblock In \emph{Uncertainty in Artificial Intelligence}, pages 493--503.
  PMLR, 2021.

\bibitem[Fern{\'a}ndez-Blanco et~al.(2021)Fern{\'a}ndez-Blanco, Morales,
  Pineda, and Porras]{fernandez2021inverse}
Ricardo Fern{\'a}ndez-Blanco, Juan~Miguel Morales, Salvador Pineda, and
  {\'A}lvaro Porras.
\newblock Inverse optimization with kernel regression: Application to the power
  forecasting and bidding of a fleet of electric vehicles.
\newblock \emph{Computers \& Operations Research}, 134:\penalty0 105405, 2021.

\bibitem[Foulds et~al.(1992)Foulds, Haugland, and
  J{\o}rnsten]{foulds1992bilinear}
Leslie~Richard Foulds, Dag Haugland, and Kurt J{\o}rnsten.
\newblock A bilinear approach to the pooling problem.
\newblock \emph{Optimization}, 24\penalty0 (1-2):\penalty0 165--180, 1992.

\bibitem[Frazier(2018)]{frazier2018tutorial}
Peter~I Frazier.
\newblock A tutorial on bayesian optimization.
\newblock \emph{arXiv preprint arXiv:1807.02811}, 2018.

\bibitem[Frazier and Wang(2016)]{frazier2016bayesian}
Peter~I Frazier and Jialei Wang.
\newblock Bayesian optimization for materials design.
\newblock \emph{Information science for materials discovery and design}, pages
  45--75, 2016.

\bibitem[Garc{\'\i}a~S{\'a}nchez and
  Torres~S{\'a}ez(2014)]{garcia2014comparison}
Carlos~Eduardo Garc{\'\i}a~S{\'a}nchez and Rodrigo~Gonzalo Torres~S{\'a}ez.
\newblock Comparison and analysis of objective functions in flux balance
  analysis.
\newblock \emph{Biotechnology progress}, 30\penalty0 (5):\penalty0 985--991,
  2014.

\bibitem[Garc{\'\i}a~S{\'a}nchez et~al.(2012)Garc{\'\i}a~S{\'a}nchez,
  Vargas~Garc{\'\i}a, and Torres~S{\'a}ez]{garcia2012predictive}
Carlos~Eduardo Garc{\'\i}a~S{\'a}nchez, C{\'e}sar~Augusto Vargas~Garc{\'\i}a,
  and Rodrigo~Gonzalo Torres~S{\'a}ez.
\newblock Predictive potential of flux balance analysis of saccharomyces
  cerevisiae using as optimization function combinations of cell compartmental
  objectives.
\newblock 2012.

\bibitem[Gardner et~al.(2018)Gardner, Pleiss, Weinberger, Bindel, and
  Wilson]{gardner2018gpytorch}
Jacob Gardner, Geoff Pleiss, Kilian~Q Weinberger, David Bindel, and Andrew~G
  Wilson.
\newblock Gpytorch: Blackbox matrix-matrix gaussian process inference with gpu
  acceleration.
\newblock \emph{Advances in neural information processing systems}, 31, 2018.

\bibitem[Ghate(2020)]{ghate2020imputing}
Archis Ghate.
\newblock Imputing radiobiological parameters of the linear-quadratic
  dose-response model from a radiotherapy fractionation plan.
\newblock \emph{Physics in Medicine \& Biology}, 65\penalty0 (22):\penalty0
  225009, 2020.

\bibitem[Ghobadi and Mahmoudzadeh(2021)]{ghobadi2021inferring}
Kimia Ghobadi and Houra Mahmoudzadeh.
\newblock Inferring linear feasible regions using inverse optimization.
\newblock \emph{European Journal of Operational Research}, 290\penalty0
  (3):\penalty0 829--843, 2021.

\bibitem[Goodfellow et~al.(2016)Goodfellow, Bengio, and
  Courville]{Goodfellow-et-al-2016}
Ian Goodfellow, Yoshua Bengio, and Aaron Courville.
\newblock \emph{Deep Learning}.
\newblock MIT Press, 2016.
\newblock \url{http://www.deeplearningbook.org}.

\bibitem[Greenhill et~al.(2020)Greenhill, Rana, Gupta, Vellanki, and
  Venkatesh]{greenhill2020bayesian}
Stewart Greenhill, Santu Rana, Sunil Gupta, Pratibha Vellanki, and Svetha
  Venkatesh.
\newblock Bayesian optimization for adaptive experimental design: A review.
\newblock \emph{IEEE access}, 8:\penalty0 13937--13948, 2020.

\bibitem[Gupta and Zhang(2022)]{gupta2022decomposition}
Rishabh Gupta and Qi~Zhang.
\newblock Decomposition and adaptive sampling for data-driven inverse linear
  optimization.
\newblock \emph{INFORMS Journal on Computing}, 34\penalty0 (5):\penalty0
  2720--2735, 2022.

\bibitem[Gupta and Zhang(2023)]{gupta2023efficient}
Rishabh Gupta and Qi~Zhang.
\newblock Efficient learning of decision-making models: A penalty block
  coordinate descent algorithm for data-driven inverse optimization.
\newblock \emph{Computers \& Chemical Engineering}, 170:\penalty0 108123, 2023.

\bibitem[{Gurobi Optimization, LLC}(2023)]{gurobi}
{Gurobi Optimization, LLC}.
\newblock {Gurobi Optimizer Reference Manual}, 2023.
\newblock URL \url{https://www.gurobi.com}.

\bibitem[Haverly(1978)]{Haverly1978studies}
Co~A Haverly.
\newblock Studies of the behavior of recursion for the pooling problem.
\newblock \emph{Acm sigmap bulletin}, \penalty0 (25):\penalty0 19--28, 1978.

\bibitem[Hutter et~al.(2011)Hutter, Hoos, and
  Leyton-Brown]{hutter2011sequential}
Frank Hutter, Holger~H Hoos, and Kevin Leyton-Brown.
\newblock Sequential model-based optimization for general algorithm
  configuration.
\newblock In \emph{Learning and Intelligent Optimization: 5th International
  Conference, LION 5, Rome, Italy, January 17-21, 2011. Selected Papers 5},
  pages 507--523. Springer, 2011.

\bibitem[Huyer and Neumaier(2008)]{huyer2008snobfit}
Waltraud Huyer and Arnold Neumaier.
\newblock Snobfit--stable noisy optimization by branch and fit.
\newblock \emph{ACM Transactions on Mathematical Software (TOMS)}, 35\penalty0
  (2):\penalty0 1--25, 2008.

\bibitem[Iyengar and Kang(2005)]{iyengar2005inverse}
Garud Iyengar and Wanmo Kang.
\newblock Inverse conic programming with applications.
\newblock \emph{Operations Research Letters}, 33\penalty0 (3):\penalty0
  319--330, 2005.

\bibitem[Jones et~al.(1998)Jones, Schonlau, and Welch]{jones1998efficient}
Donald~R Jones, Matthias Schonlau, and William~J Welch.
\newblock Efficient global optimization of expensive black-box functions.
\newblock \emph{Journal of Global optimization}, 13:\penalty0 455--492, 1998.

\bibitem[Joshi et~al.(2006)Joshi, Seidel-Morgenstern, and
  Kremling]{joshi2006exploiting}
Milind Joshi, Andreas Seidel-Morgenstern, and Andreas Kremling.
\newblock Exploiting the bootstrap method for quantifying parameter confidence
  intervals in dynamical systems.
\newblock \emph{Metabolic engineering}, 8\penalty0 (5):\penalty0 447--455,
  2006.

\bibitem[Keshavarz et~al.(2011)Keshavarz, Wang, and
  Boyd]{keshavarz2011imputing}
Arezou Keshavarz, Yang Wang, and Stephen Boyd.
\newblock Imputing a convex objective function.
\newblock In \emph{2011 IEEE international symposium on intelligent control},
  pages 613--619. IEEE, 2011.

\bibitem[Kieffer et~al.(2017)Kieffer, Danoy, Bouvry, and
  Nagih]{kieffer2017bayesian}
Emmanuel Kieffer, Gr{\'e}goire Danoy, Pascal Bouvry, and Anass Nagih.
\newblock Bayesian optimization approach of general bi-level problems.
\newblock In \emph{Proceedings of the Genetic and Evolutionary Computation
  Conference Companion}, pages 1614--1621, 2017.

\bibitem[King et~al.(2016)King, Lu, Dr{\"a}ger, Miller, Federowicz, Lerman,
  Ebrahim, Palsson, and Lewis]{king2016bigg}
Zachary~A King, Justin Lu, Andreas Dr{\"a}ger, Philip Miller, Stephen
  Federowicz, Joshua~A Lerman, Ali Ebrahim, Bernhard~O Palsson, and Nathan~E
  Lewis.
\newblock Bigg models: A platform for integrating, standardizing and sharing
  genome-scale models.
\newblock \emph{Nucleic acids research}, 44\penalty0 (D1):\penalty0 D515--D522,
  2016.

\bibitem[Kreutz and Timmer(2009)]{kreutz2009systems}
Clemens Kreutz and Jens Timmer.
\newblock Systems biology: experimental design.
\newblock \emph{The FEBS journal}, 276\penalty0 (4):\penalty0 923--942, 2009.

\bibitem[Kudva et~al.(2022)Kudva, Sorourifar, and
  Paulson]{kudva2022constrained}
Akshay Kudva, Farshud Sorourifar, and Joel~A Paulson.
\newblock Constrained robust bayesian optimization of expensive noisy black-box
  functions with guaranteed regret bounds.
\newblock \emph{AIChE Journal}, 68\penalty0 (12):\penalty0 e17857, 2022.

\bibitem[Kudva et~al.(2024)Kudva, Tang, and Paulson]{kudva2024robust}
Akshay Kudva, Wei-Ting Tang, and Joel~A Paulson.
\newblock Robust bayesian optimization for flexibility analysis of expensive
  simulation-based models with rigorous uncertainty bounds.
\newblock \emph{Computers \& Chemical Engineering}, 181:\penalty0 108515, 2024.

\bibitem[Lee et~al.(2009)Lee, Burd, Liu, Almaas, Wiest, Barab{\'a}si, Oltvai,
  and Kapatral]{lee2009comparative}
Deok-Sun Lee, Henry Burd, Jiangxia Liu, Eivind Almaas, Olaf Wiest,
  Albert-L{\'a}szl{\'o} Barab{\'a}si, Zolt{\'a}n~N Oltvai, and Vinayak
  Kapatral.
\newblock Comparative genome-scale metabolic reconstruction and flux balance
  analysis of multiple staphylococcus aureus genomes identify novel
  antimicrobial drug targets.
\newblock \emph{Journal of bacteriology}, 191\penalty0 (12):\penalty0
  4015--4024, 2009.

\bibitem[Lee et~al.(2004)Lee, Hwang, Yokoyama, Stephanopoulos, Stephanopoulos,
  and Yarmush]{lee2004identification}
Kyongbum Lee, Daehee Hwang, Tadaaki Yokoyama, George Stephanopoulos, Gregory~N
  Stephanopoulos, and Martin~L Yarmush.
\newblock Identification of optimal classification functions for biological
  sample and state discrimination from metabolic profiling data.
\newblock \emph{Bioinformatics}, 20\penalty0 (6):\penalty0 959--969, 2004.

\bibitem[Lee and Grossmann(2003)]{lee2003global}
Sangbum Lee and Ignacio~E Grossmann.
\newblock Global optimization of nonlinear generalized disjunctive programming
  with bilinear equality constraints: applications to process networks.
\newblock \emph{Computers \& chemical engineering}, 27\penalty0 (11):\penalty0
  1557--1575, 2003.

\bibitem[Lewis et~al.(2010)Lewis, Hixson, Conrad, Lerman, Charusanti,
  Polpitiya, Adkins, Schramm, Purvine, Lopez-Ferrer, et~al.]{lewis2010omic}
Nathan~E Lewis, Kim~K Hixson, Tom~M Conrad, Joshua~A Lerman, Pep Charusanti,
  Ashoka~D Polpitiya, Joshua~N Adkins, Gunnar Schramm, Samuel~O Purvine, Daniel
  Lopez-Ferrer, et~al.
\newblock Omic data from evolved e. coli are consistent with computed optimal
  growth from genome-scale models.
\newblock \emph{Molecular systems biology}, 6\penalty0 (1):\penalty0 390, 2010.

\bibitem[Li(2021)]{li2021inverse}
Jonathan Yu-Meng Li.
\newblock Inverse optimization of convex risk functions.
\newblock \emph{Management Science}, 67\penalty0 (11):\penalty0 7113--7141,
  2021.

\bibitem[Maiwald et~al.(2016)Maiwald, Hass, Steiert, Vanlier, Engesser, Raue,
  Kipkeew, Bock, Kaschek, Kreutz, et~al.]{maiwald2016driving}
Tim Maiwald, Helge Hass, Bernhard Steiert, Joep Vanlier, Raphael Engesser,
  Andreas Raue, Friederike Kipkeew, Hans~H Bock, Daniel Kaschek, Clemens
  Kreutz, et~al.
\newblock Driving the model to its limit: profile likelihood based model
  reduction.
\newblock \emph{PloS one}, 11\penalty0 (9):\penalty0 e0162366, 2016.

\bibitem[Misener and Floudas(2009)]{misener2009advances}
Ruth Misener and Christodoulos~A Floudas.
\newblock Advances for the pooling problem: Modeling, global optimization, and
  computational studies.
\newblock \emph{Applied and Computational Mathematics}, 8\penalty0
  (1):\penalty0 3--22, 2009.

\bibitem[Mockus(1994)]{mockus1994application}
Jonas Mockus.
\newblock Application of bayesian approach to numerical methods of global and
  stochastic optimization.
\newblock \emph{Journal of Global Optimization}, 4:\penalty0 347--365, 1994.

\bibitem[Mohajerin~Esfahani et~al.(2018)Mohajerin~Esfahani,
  Shafieezadeh-Abadeh, Hanasusanto, and Kuhn]{mohajerin2018data}
Peyman Mohajerin~Esfahani, Soroosh Shafieezadeh-Abadeh, Grani~A Hanasusanto,
  and Daniel Kuhn.
\newblock Data-driven inverse optimization with imperfect information.
\newblock \emph{Mathematical Programming}, 167:\penalty0 191--234, 2018.

\bibitem[Nagrath et~al.(2010)Nagrath, Avila-Elchiver, Berthiaume, Tilles,
  Messac, and Yarmush]{nagrath2010soft}
Deepak Nagrath, Marco Avila-Elchiver, Fran{\c{c}}ois Berthiaume, Arno~W Tilles,
  Achille Messac, and Martin~L Yarmush.
\newblock Soft constraints-based multiobjective framework for flux balance
  analysis.
\newblock \emph{Metabolic engineering}, 12\penalty0 (5):\penalty0 429--445,
  2010.

\bibitem[Neale and Miller(1997)]{neale1997use}
Michael~C Neale and Michael~B Miller.
\newblock The use of likelihood-based confidence intervals in genetic models.
\newblock \emph{Behavior genetics}, 27:\penalty0 113--120, 1997.

\bibitem[Nilsson and Nielsen(2017)]{nilsson2017genome}
Avlant Nilsson and Jens Nielsen.
\newblock Genome scale metabolic modeling of cancer.
\newblock \emph{Metabolic engineering}, 43:\penalty0 103--112, 2017.

\bibitem[Orth et~al.(2010{\natexlab{a}})Orth, Fleming, and
  Palsson]{orth2010reconstruction}
Jeffrey~D Orth, Ronan~MT Fleming, and Bernhard~{\O} Palsson.
\newblock Reconstruction and use of microbial metabolic networks: the core
  escherichia coli metabolic model as an educational guide.
\newblock \emph{EcoSal plus}, 4\penalty0 (1):\penalty0 10--1128,
  2010{\natexlab{a}}.

\bibitem[Orth et~al.(2010{\natexlab{b}})Orth, Thiele, and
  Palsson]{orth2010flux}
Jeffrey~D Orth, Ines Thiele, and Bernhard~{\O} Palsson.
\newblock What is flux balance analysis?
\newblock \emph{Nature biotechnology}, 28\penalty0 (3):\penalty0 245--248,
  2010{\natexlab{b}}.

\bibitem[Paciorek and Schervish(2003)]{paciorek2003nonstationary}
Christopher Paciorek and Mark Schervish.
\newblock Nonstationary covariance functions for gaussian process regression.
\newblock \emph{Advances in neural information processing systems}, 16, 2003.

\bibitem[Pardelha et~al.(2012)Pardelha, Albuquerque, Reis, Dias, and
  Oliveira]{pardelha2012flux}
Filipa Pardelha, Maria~GE Albuquerque, Maria~AM Reis, Jo{\~a}o~ML Dias, and Rui
  Oliveira.
\newblock Flux balance analysis of mixed microbial cultures: Application to the
  production of polyhydroxyalkanoates from complex mixtures of volatile fatty
  acids.
\newblock \emph{Journal of biotechnology}, 162\penalty0 (2-3):\penalty0
  336--345, 2012.

\bibitem[Raue et~al.(2009)Raue, Kreutz, Maiwald, Bachmann, Schilling,
  Klingm{\"u}ller, and Timmer]{raue2009structural}
Andreas Raue, Clemens Kreutz, Thomas Maiwald, Julie Bachmann, Marcel Schilling,
  Ursula Klingm{\"u}ller, and Jens Timmer.
\newblock Structural and practical identifiability analysis of partially
  observed dynamical models by exploiting the profile likelihood.
\newblock \emph{Bioinformatics}, 25\penalty0 (15):\penalty0 1923--1929, 2009.

\bibitem[Rios and Sahinidis(2013)]{rios2013derivative}
Luis~Miguel Rios and Nikolaos~V Sahinidis.
\newblock Derivative-free optimization: a review of algorithms and comparison
  of software implementations.
\newblock \emph{Journal of Global Optimization}, 56\penalty0 (3):\penalty0
  1247--1293, 2013.

\bibitem[Savinell and Palsson(1992)]{savinell1992network}
Joanne~M Savinell and Bernhard~O Palsson.
\newblock Network analysis of intermediary metabolism using linear
  optimization. i. development of mathematical formalism.
\newblock \emph{Journal of theoretical biology}, 154\penalty0 (4):\penalty0
  421--454, 1992.

\bibitem[Schaefer(2009)]{schaefer2009inverse}
Andrew~J Schaefer.
\newblock Inverse integer programming.
\newblock \emph{Optimization Letters}, 3:\penalty0 483--489, 2009.

\bibitem[Schoemaker(1991)]{schoemaker1991quest}
Paul~JH Schoemaker.
\newblock The quest for optimality: A positive heuristic of science?
\newblock \emph{Behavioral and brain sciences}, 14\penalty0 (2):\penalty0
  205--215, 1991.

\bibitem[Schuetz et~al.(2007)Schuetz, Kuepfer, and
  Sauer]{schuetz2007systematic}
Robert Schuetz, Lars Kuepfer, and Uwe Sauer.
\newblock Systematic evaluation of objective functions for predicting
  intracellular fluxes in escherichia coli.
\newblock \emph{Molecular systems biology}, 3\penalty0 (1):\penalty0 119, 2007.

\bibitem[Schuetz et~al.(2012)Schuetz, Zamboni, Zampieri, Heinemann, and
  Sauer]{schuetz2012multidimensional}
Robert Schuetz, Nicola Zamboni, Mattia Zampieri, Matthias Heinemann, and Uwe
  Sauer.
\newblock Multidimensional optimality of microbial metabolism.
\newblock \emph{Science}, 336\penalty0 (6081):\penalty0 601--604, 2012.

\bibitem[Shahriari et~al.(2015)Shahriari, Swersky, Wang, Adams, and
  De~Freitas]{shahriari2015taking}
Bobak Shahriari, Kevin Swersky, Ziyu Wang, Ryan~P Adams, and Nando De~Freitas.
\newblock Taking the human out of the loop: A review of bayesian optimization.
\newblock \emph{Proceedings of the IEEE}, 104\penalty0 (1):\penalty0 148--175,
  2015.

\bibitem[Sinha et~al.(2017)Sinha, Malo, and Deb]{sinha2017review}
Ankur Sinha, Pekka Malo, and Kalyanmoy Deb.
\newblock A review on bilevel optimization: From classical to evolutionary
  approaches and applications.
\newblock \emph{IEEE transactions on evolutionary computation}, 22\penalty0
  (2):\penalty0 276--295, 2017.

\bibitem[Snoek et~al.(2012)Snoek, Larochelle, and Adams]{snoek2012practical}
Jasper Snoek, Hugo Larochelle, and Ryan~P Adams.
\newblock Practical bayesian optimization of machine learning algorithms.
\newblock \emph{Advances in neural information processing systems}, 25, 2012.

\bibitem[Snoek et~al.(2015)Snoek, Rippel, Swersky, Kiros, Satish, Sundaram,
  Patwary, Prabhat, and Adams]{snoek2015scalable}
Jasper Snoek, Oren Rippel, Kevin Swersky, Ryan Kiros, Nadathur Satish,
  Narayanan Sundaram, Mostofa Patwary, Mr~Prabhat, and Ryan Adams.
\newblock Scalable bayesian optimization using deep neural networks.
\newblock In \emph{International conference on machine learning}, pages
  2171--2180. PMLR, 2015.

\bibitem[Sobel(1982)]{sobel1982asymptotic}
Michael~E Sobel.
\newblock Asymptotic confidence intervals for indirect effects in structural
  equation models.
\newblock \emph{Sociological methodology}, 13:\penalty0 290--312, 1982.

\bibitem[Springenberg et~al.(2016)Springenberg, Klein, Falkner, and
  Hutter]{springenberg2016bayesian}
Jost~Tobias Springenberg, Aaron Klein, Stefan Falkner, and Frank Hutter.
\newblock Bayesian optimization with robust bayesian neural networks.
\newblock \emph{Advances in neural information processing systems}, 29, 2016.

\bibitem[Srinivas et~al.(2009)Srinivas, Krause, Kakade, and
  Seeger]{srinivas2009gaussian}
Niranjan Srinivas, Andreas Krause, Sham~M Kakade, and Matthias Seeger.
\newblock Gaussian process optimization in the bandit setting: No regret and
  experimental design.
\newblock \emph{arXiv preprint arXiv:0912.3995}, 2009.

\bibitem[Talbi(2013)]{talbi2013taxonomy}
El-Ghazali Talbi.
\newblock A taxonomy of metaheuristics for bi-level optimization.
\newblock In \emph{Metaheuristics for bi-level optimization}, pages 1--39.
  Springer, 2013.

\bibitem[Uygun et~al.(2007)Uygun, Matthew, and Huang]{uygun2007investigation}
Korkut Uygun, Howard~WT Matthew, and Yinlun Huang.
\newblock Investigation of metabolic objectives in cultured hepatocytes.
\newblock \emph{Biotechnology and bioengineering}, 97\penalty0 (3):\penalty0
  622--637, 2007.

\bibitem[Wang(2009)]{wang2009cutting}
Lizhi Wang.
\newblock Cutting plane algorithms for the inverse mixed integer linear
  programming problem.
\newblock \emph{Operations research letters}, 37\penalty0 (2):\penalty0
  114--116, 2009.

\bibitem[Wieland et~al.(2021)Wieland, Hauber, Rosenblatt, T{\"o}nsing, and
  Timmer]{wieland2021structural}
Franz-Georg Wieland, Adrian~L Hauber, Marcus Rosenblatt, Christian T{\"o}nsing,
  and Jens Timmer.
\newblock On structural and practical identifiability.
\newblock \emph{Current Opinion in Systems Biology}, 25:\penalty0 60--69, 2021.

\bibitem[Williams and Rasmussen(2006)]{williams2006gaussian}
Christopher~KI Williams and Carl~Edward Rasmussen.
\newblock \emph{Gaussian processes for machine learning}, volume~2.
\newblock MIT press Cambridge, MA, 2006.

\bibitem[Yu et~al.(2023)Yu, Wang, and Dong]{yu2023learning}
Shi Yu, Haoran Wang, and Chaosheng Dong.
\newblock Learning risk preferences from investment portfolios using inverse
  optimization.
\newblock \emph{Research in International Business and Finance}, 64:\penalty0
  101879, 2023.

\bibitem[Zhang and Xu(2010)]{zhang2010inverse}
Jianzhong Zhang and Chengxian Xu.
\newblock Inverse optimization for linearly constrained convex separable
  programming problems.
\newblock \emph{European Journal of Operational Research}, 200\penalty0
  (3):\penalty0 671--679, 2010.

\bibitem[Zhang and Zhang(2010)]{zhang2010augmented}
Jianzhong Zhang and Liwei Zhang.
\newblock An augmented lagrangian method for a class of inverse quadratic
  programming problems.
\newblock \emph{Applied Mathematics and Optimization}, 61\penalty0
  (1):\penalty0 57, 2010.

\bibitem[Zhang et~al.(2018)Zhang, Pourazarm, Cassandras, and
  Paschalidis]{zhang2018price}
Jing Zhang, Sepideh Pourazarm, Christos~G Cassandras, and Ioannis~Ch
  Paschalidis.
\newblock The price of anarchy in transportation networks: Data-driven
  evaluation and reduction strategies.
\newblock \emph{Proceedings of the IEEE}, 106\penalty0 (4):\penalty0 538--553,
  2018.

\bibitem[Zhao et~al.(2016)Zhao, Stettner, Reznik, Paschalidis, and
  Segr{\`e}]{zhao2016mapping}
Qi~Zhao, Arion~I Stettner, Ed~Reznik, Ioannis~Ch Paschalidis, and Daniel
  Segr{\`e}.
\newblock Mapping the landscape of metabolic goals of a cell.
\newblock \emph{Genome biology}, 17:\penalty0 1--11, 2016.

\end{thebibliography}

\end{document}